\documentclass[12pt,letterpaper]{amsart}
\usepackage{amscd,amsxtra}
\usepackage[all]{xy}

\newtheorem{theorem}{Theorem}[section]
\newtheorem{corollary}[theorem]{Corollary}
\newtheorem{lemma}[theorem]{Lemma}
\newtheorem{proposition}[theorem]{Proposition}

\theoremstyle{definition}
\newtheorem{definition}[theorem]{Definition}
\newtheorem{remark}[theorem]{Remark}
\newtheorem{example}[theorem]{Example}

\DeclareMathOperator{\Coh}{\mathsf{Coh}} 
\DeclareMathOperator{\Ext}{\mathsf{Ext}}
\DeclareMathOperator{\Hom}{\mathsf{Hom}}
\DeclareMathOperator{\shom}{\mathsf{hom}}
\DeclareMathOperator{\cone}{\mathsf{cone}}
\DeclareMathOperator{\Pic}{Pic}
\DeclareMathOperator{\rk}{rk}
\DeclareMathOperator{\coker}{coker}
\DeclareMathOperator{\Hilb}{\mathsf{Hilb}}
\DeclareMathOperator{\Quot}{\mathsf{Quot}} 
 
\DeclareMathOperator{\supp}{supp} 
\DeclareMathOperator{\depth}{depth} 
 
\DeclareMathOperator{\cha}{char}
\DeclareMathOperator{\Id}{\mathsf{Id}}
\DeclareMathOperator{\SL}{\mathsf{SL}}
\newcommand{\dtens}{\stackrel{\boldsymbol{L}}{\otimes}} 
\newcommand{\local}[2]{\boldsymbol{k}[[#1,#2]]/(#1\cdot #2)}

\begin{document}
\title
[Fourier-Mukai transforms and semi-stable sheaves]
{Fourier-Mukai transforms and semi-stable sheaves on nodal Weierstra{\ss}
cubics} 

\author{Igor Burban}
\address{%
Institut de Math\'ematiques de Jussieu, 
Universit\'e  Pierre et Marie Curie -- Paris VI,
4 Place Jussieu, 75252 PARIS Cedex 05, France
}
\email{burban@math.jussieu.fr}

\author{Bernd Kreu{\ss}ler}
\address{%
Mary Immaculate College, South Circular Road, Limerick, Ireland
}
\email{bernd.kreussler@mic.ul.ie}

\thanks{%
Both authors were partially supported by DFG Schwerpunkt 
``Globale Methoden in der komplexen Geometrie''.}


\dedicatory{Dedicated to Gert--Martin Greuel on his sixtieth birthday}

\begin{abstract}
We completely describe all semi-stable torsion free sheaves of
degree zero on nodal cubic curves using the technique of Fourier-Mukai
transforms. The Fourier-Mukai images of such sheaves are torsion sheaves of
finite length, which we compute explicitly. We show that the twist functors,
which are associated to the structure sheaf $\mathcal{O}$ and the structure
sheaf $\boldsymbol{k}(p_{0})$ of a smooth point $p_{0}$, generate an
$\SL(2,\mathbb{Z})$-action (up to shifts) on the bounded derived category of
coherent sheaves on any Weiersta{\ss} cubic.
\end{abstract}

\maketitle

\section{Introduction}
\label{sec:intro}

In recent years, derived categories of coherent sheaves on smooth projective
varieties and their groups of auto-equivalences attracted a lot of
interest. This was mainly driven by Kontsevich's homological mirror
symmetry conjecture \cite{Kontsevich}. It turns out \cite{BondalOrlov} that
interesting auto-equivalences can exist only if the canonical sheaf and its
dual are not ample. From a mirror symmetry perspective, the most interesting
case is the case with trivial canonical sheaf. One particular class of
varieties with trivial canonical bundle are the abelian varieties. If $X$ is
an abelian variety, $X^{\vee}$ its dual, $\mathcal{P}$ a Poincar\'e bundle on
$X\times X^{\vee}$ and $\pi_{1}, \pi_{2}$ are the two projections, Mukai
\cite{Mukai} has shown that the functor 
$$\Phi_{\mathcal{P}}:D^{b}(X) \rightarrow D^{b}(X^{\vee}),
 \quad \Phi_{\mathcal{P}}(F) = 
       \boldsymbol{R}\pi_{2\ast}(\mathcal{P}\dtens \pi^{\ast}_{1}F)$$
is an exact equivalence of categories. Nowadays, such functors are called
\textsl{Fourier-Mukai transforms}. If $X$ is a smooth elliptic curve (or
more generally: a principally polarised abelian variety), one has an
isomorphism $X\cong X^{\vee}$. In this case, the above functor induces an
auto-equivalence of the derived category of $X$. In general, a result of
Orlov \cite{Orlov} says that any equivalence between bounded derived
categories of coherent sheaves on smooth projective varieties $X$ and $Y$ is a
Fourier-Mukai transform, where we allow $\mathcal{P}$ to be replaced by an
arbitrary object of $D^{b}(X\times Y)$.

On the other hand, not so much is known in the singular case. In this paper,
we study the first non-trivial example: singular Weierstra{\ss} cubics. 
Our principal tool is a Fourier-Mukai transform
$\mathbb{F}:D^{b}(\boldsymbol{E})\rightarrow D^{b}(\boldsymbol{E})$, which we
define and study with techniques introduced by Seidel and Thomas
\cite{SeidelThomas} (see also \cite{GorodentsevRudakov}, \cite{LenzMeltzer}
and \cite{Meltzer}). Our results are of interest in the study of relative
Fourier-Mukai transforms on elliptic fibrations, see for example \cite{BBRP}
and \cite{BridgelandElliptic}.
Interesting applications of Fourier-Mukai transforms on Weierstra{\ss} cubics
to Calogero-Moser systems can be found in \cite{BZN}. 

In a series of papers, Friedman, Morgan and Witten \cite{FMW}, \cite{FMmin},
\cite{FM} studied, among other things, semi-stable vector bundles on
Weierstra{\ss} cubics. We extend their results and give a complete and
explicit description of all semi-stable torsion free sheaves of degree zero on
a nodal Weierstra{\ss} cubic $\boldsymbol{E}$. A description of all stable
vector bundles (of any degree) on $\boldsymbol{E}$ can be found in \cite{UMZ}.
Although a description of all torsion free sheaves on
$\boldsymbol{E}$ is available, it is a non-trivial problem to find all
semi-stable torsion free sheaves among them. With the aid of the Fourier-Mukai
transform $\mathbb{F}$ we are able to translate this problem into a problem
which is known as a matrix problem of Gelfand type in representation theory. 
We use the solution of this problem which was given in \cite{GelfandPonomarev}
and provide an explicit description of the functor $\mathbb{F}$ on the sheaves
in question. 
An explicit description of $\mathbb{F}$ on semi-stable torsion free sheaves
seems to be worthwhile, because the interplay between such sheaves and torsion
sheaves was used as an efficient technical tool by Friedman, Morgan and
others. As a result, we obtain a clear description of all semi-stable torsion
free sheaves of degree zero and their Fourier-Mukai images. 

It is instructive to associate to such a sheaf (or their Fourier-Mukai image)
the so-called band and string diagram from representation theory. 
They are particularly useful in the study of the Fourier-Mukai image of the
dual of a sheaf. This leads us to Matlis duality over local Gorenstein
rings. In the nodal case, Matlis duality is described by reversing the
direction of all arrows in the band and string diagrams. 

The structure and more detailed content of this paper is the following. 
In section \ref{sec:fm} we use twist functors which are defined by spherical
objects to give a definition of the Fourier-Mukai functor
$\mathbb{F}:D^{b}(\boldsymbol{E})\rightarrow D^{b}(\boldsymbol{E})$. 
One of the main results in this section is (theorem \ref{main}):
$$\mathbb{F}\circ \mathbb{F} = i^{\ast}[1],$$ where
$i:\boldsymbol{E}\rightarrow\boldsymbol{E}$ is the involution which is induced
by taking the inverse in the group structure on the regular locus of
$\boldsymbol{E}$. 

This theorem and a braid group relation between two particular twist functors
allow us to define an action (up to translation and in a weak sense) of
$\SL(2,\mathbb{Z})$  on $D^{b}(\boldsymbol{E})$. 
If we let auto-equivalences of $D^{b}(\boldsymbol{E})$ act on pairs formed by 
rank and degree, which can be defined for complexes as well, we obtain a
surjection from the group of exact auto-equivalences onto $\SL(2,\mathbb{Z})$. 
This is parallel to the smooth case. 
But, in contrast, we do not know whether the whole group of
auto-equivalences of $D^{b}(\boldsymbol{E})$ is generated by 
this $\SL(2,\mathbb{Z})$-action and automorphisms of $\boldsymbol{E}$. 
Clearly, all shifts and twists by line bundles are obtained this way. 

If $\boldsymbol{E}$ is a reducible curve of arithmetic genus one, for example
a cycle of rational curves, most of the techniques which are used so far are
applicable as well. However, in this case it is more complicated to describe
$\mathbb{F}\circ \mathbb{F}$. Furthermore, if $\boldsymbol{E}$ has $n>1$
components, the rank lives in $\mathbb{Z}^{n}$, because it needs to be
constant on the irreducible components only. 
Hence, the auto-equivalences of $D^{b}(\boldsymbol{E})$
naturally act on $\mathbb{Z}^{n+1}$ and we expect to be forced to replace
$\SL(2,\mathbb{Z})$ by another group.

The second main result (theorem \ref{ss}) says that
$\mathbb{F}$ induces an equivalence of categories: 
$$ \left\{
  \begin{array}{c}
    \text{semi-stable torsion free}\\ 
    \text{sheaves of degree zero}
  \end{array}         
  \right\}  
  \longrightarrow
  \left\{
  \begin{array}{c}
    \text{coherent torsion sheaves}
  \end{array}         
\right\}.
$$
On this category, our functor $\mathbb{F}$ coincides with the functor
introduced by Teodorescu \cite{Teodorescu}, see also \cite{FMmin}. 

In section \ref{sec:tf} we describe vector bundles on $\boldsymbol{E}$ as
direct images of line bundles under cyclic Galois covers, provided
$\boldsymbol{E}$ is a nodal Weierstra{\ss} cubic. The results are
similar to Atiyah's classical results for the smooth case \cite{Atiyah}. In
addition, we describe torsion free sheaves which are not locally free as
direct images of line bundles via a finite morphism from chains of rational
curves. In the smooth case, such sheaves do not exist. It can be derived
from \cite{DrozdGreuel} that the sheaves described 
this way comprise all torsion free sheaves on $\boldsymbol{E}$. To do so, one
has to relate their description of sheaves via gluing of bundles to the
description by direct images of line bundles. This is possible in
characteristic zero only. To prove results which are similar to ours over a
field of arbitrary characteristic, one had to avoid the use of this description
by direct images. 

In section \ref{sec:tors} we recall the results of Gelfand, Ponomarev
\cite{GelfandPonomarev} and Nazarova, Roiter \cite{NazRoi} on the
classification of $R$-modules of finite length. In 
contrast to the case of a regular local ring of dimension one, where
indecomposable modules of finite length are determined by their length,
over the complete local Gorenstein ring $R=\local{x}{y}$ there are plenty of
such modules. 
As a useful tool to describe them, we introduce the band and string diagrams,
which were already used in \cite{GelfandPonomarev}.

The main result of this paper (theorem \ref{description}) is proven in section
\ref{sec:semistable}. We explicitly compute the Fourier-Mukai images of
certain torsion free sheaves. Using the results of sections \ref{sec:fm} and
\ref{sec:tors} we are able to conclude that these sheaves are precisely the
semi-stable torsion free sheaves of degree zero on a nodal Weierstra{\ss}
cubic $\boldsymbol{E}$. 

Finally, in section \ref{sec:matlis}, we study the relationship between
$\mathbb{F}(E)$ and $\mathbb{F}(E^{\vee})$, where $E$ is a semi-stable torsion
free sheaf of degree zero on $\boldsymbol{E}$. After briefly recalling Matlis
duality we show (theorem \ref{dual}): 
$$\mathbb{F}(E^{\vee}) \cong \mathbb{F}(E)^{\star},$$
where $M^{\star}$ denotes the twisted Matlis dual (definition \ref{twistedmd})
of a module $M$ of finite length. 

Throughout this paper, $\boldsymbol{k}$ denotes an algebraically closed field
of characteristic zero. All schemes are schemes over $\boldsymbol{k}$.

\textbf{Acknowledgement.} Both authors like to thank G.-M.\/ Greuel for his
constant support and Yu.\/ Drozd, H.\/ Krause, H.\/ Lenzing and R.\/ Rouquier
for helpful discussions about results of this paper.

\section{Fourier-Mukai functors on Weierstra{\ss} cubics}
\label{sec:fm}

For abelian varieties $X$, Mukai \cite{Mukai}, Thm. 2.2, has shown that the
integral transform  
$\Phi_{\mathcal{P}}:D^{b}(X) \rightarrow D^{b}(X^{\vee})$ which is given by 
$$\Phi_{\mathcal{P}}(F) :=
\boldsymbol{R}\pi_{2\ast}(\mathcal{P}\dtens \pi^{\ast}_{1}F)$$
is an equivalence of categories. Here, $\mathcal{P}\in\Coh(X\times X^{\vee})$
is a Poincar\'e bundle and $X^{\vee}$ is the dual abelian variety. Nowadays,
such an integral transform with arbitrary $\mathcal{P}\in D^{b}(X\times Y)$ is
called a \emph{Fourier-Mukai transform} if it is an equivalence of categories. 

We are dealing here with the one dimensional case (elliptic curves) and
attempt to generalise this result to singular curves. In doing so, one
encounters difficulties which are caused by the fact that a Poincar\'e sheaf
$\mathcal{P}$ is not locally free if the curve is singular. To circumvent such
problems, we use techniques and results from P. Seidel and R. Thomas
\cite{SeidelThomas}. They study twist functors defined by spherical objects
and relate them to Fourier-Mukai transforms. 
These twist functors resemble the mutations which were used in the study of
exceptional vector bundles (see \cite{Drezet}, \cite{GorodentsevRudakov}).
In a more algebraic framework such functors were studied by
H.~Meltzer in \cite{Meltzer}, where they are called tubular mutations (see
also \cite{LenzMeltzer}).

\begin{definition}\label{weierstrass}
  A \emph{Weierstra{\ss} cubic} $\boldsymbol{E}$ is a plane cubic curve which
  is given in $\mathbb{P}^{2}$ by an equation
  $$y^{2}z = 4x^{3} -g_{2}xz^{2} - g_{3}z^{3},$$
  where $(x:y:z)$ are homogeneous coordinates on $\mathbb{P}^{2}$ and $g_{2},
  g_{3}\in\boldsymbol{k}$ are constants.
\end{definition}

\begin{remark}
  Weierstra{\ss} cubics are precisely the reduced and irreducible curves of
  arithmetic genus one.
  A Weierstra{\ss} cubic $\boldsymbol{E}$ has at most one singular point. It is
  singular if and only if $g_{2}^{3}=27g_{3}^{2}$. Unless $g_{2}=g_{3}=0$, the
  singularity is a node (ordinary double point), whereas in the case
  $g_{2}=g_{3}=0$ the singularity is a cusp. 
\end{remark}

If $\boldsymbol{E}$ is singular, its normalisation is a
morphism $n:\mathbb{P}^{1}\rightarrow\boldsymbol{E}$ and the canonical
morphism of sheaves $\mathcal{O}\rightarrow n_{\ast}n^{\ast}\mathcal{O}$ on
$\boldsymbol{E}$ is injective with cokernel $\boldsymbol{k}(s)$, the
structure sheaf of the singular point $s\in\boldsymbol{E}$.

Because any Weierstra{\ss} cubic $\boldsymbol{E}$ is given by a single
equation in $\mathbb{P}^{2}$, it is Gorenstein. Its dualising sheaf is
$\omega_{\boldsymbol{E}} \cong \mathcal{O}_{\boldsymbol{E}}$, because the
equation has degree three.
In particular, Serre duality \cite{Hartshorne}, II.7, holds, this means 
for any integer $i\in\mathbb{Z}$ there exists a functorial isomorphism:
$$\Ext^{i}(\mathcal{E},\mathcal{F}) \cong
\Ext^{1-i}(\mathcal{F},\mathcal{E})^{\vee},$$
where  $\mathcal{F}\in\Coh(\boldsymbol{E})$ is a coherent sheaf and 
$\mathcal{E}$  is locally free of finite rank.

We use $D^{b}(\boldsymbol{E})$ to denote the bounded derived category of
coherent sheaves on $\boldsymbol{E}$. Complexes $F^{\bullet}$ are considered to
be cochain complexes: $\cdots\rightarrow F^{-2} \rightarrow F^{-1} \rightarrow 
F^{0} \rightarrow F^{1} \rightarrow F^{2} \rightarrow \cdots$ and the shift
functor is $(F^{\bullet}[1])^{k} = F^{k+1}$. In particular, if $E$ is a
coherent sheaf, $E[-k]$ denotes the complex with $E$ at position $k$ and zero
elsewhere. 

Using induction on the length of the complex and the standard notation 
$\Hom^{i}(E^{\bullet},F^{\bullet}) = \Hom_{D^{b}(\boldsymbol{E})}(E^{\bullet},
F^{\bullet}[i])$, we deduce the following from Serre-duality in its formulation
above: for any $i\in\mathbb{Z}$ and $F\in D^{b}(\boldsymbol{E})$ and any
bounded complex of locally free sheaves $E$ there exist functorial isomorphisms
$$\Hom^{i}(E,F) \cong \Hom^{1-i}(F,E)^{\vee}.$$

On a smooth curve this is true for arbitrary $E\in
D^{b}(\boldsymbol{E})$, because any coherent sheaf on a smooth projective
variety has a finite resolution by locally free sheaves. If $\boldsymbol{E}$
is singular, this is no longer true. 

The definition of spherical objects which was given by P.\/ Seidel and R.\/
Thomas \cite{SeidelThomas}, Definition 2.14, reads in our context as follows:

\begin{definition}
  An object $E\in D^{b}(\boldsymbol{E})$ is called \emph{spherical}, if the
  following conditions are satisfied:
  \begin{itemize}
  \item[(S1)] $E$ has a finite resolution by injective quasi-coherent sheaves; 
  \item[(S2)] for any $F\in D^{b}(\boldsymbol{E})$ the total morphism spaces
    $\Hom^{\ast}(E,F)$ and $\Hom^{\ast}(F,E)$ are of finite dimension;
  \item[(S3)] 
    $$\Hom^{i}(E,E) \cong \begin{cases}
      \boldsymbol{k}&\text{ if } i=0,1\\
      0&\text{ otherwise;} 
    \end{cases}$$
  \item[(S4)] the composition map $$\Hom^{i}(F,E) \times \Hom^{1-i}(E,F)
    \rightarrow \Hom^{1}(E,E) \cong \boldsymbol{k}$$ is a non-degenerate
    pairing for any $i\in\mathbb{Z}$ and any $F\in D^{b}(\boldsymbol{E})$.
  \end{itemize}
\end{definition}

\begin{proposition}
  If $\boldsymbol{E}$ is a Weierstra{\ss} cubic and $E\in
  D^{b}(\boldsymbol{E})$ is isomorphic to a bounded complex of locally free
  sheaves, then: 
  $$E \text{ is spherical}\quad \iff \quad \Hom^{i}(E,E)\cong
  \begin{cases}
    \boldsymbol{k}&\text{ if } i=0\\
    0&\text{ if } i<0
  \end{cases}$$
\end{proposition}

\begin{proof}
  Serre duality gives (S4). The assumption and (S4) imply (S3). Now, (S2)
  follows from (S4), $\Hom^{\ast}(E,F)\cong \mathbb{H}^{\ast}(E^{\vee}\otimes
  F)$ and the standard finiteness theorem on projective schemes. Finally, (S1)
  is true because $\boldsymbol{E}$ is Gorenstein and hence any locally free
  sheaf has finite injective dimension
  (see e.g.\/ \cite{BrunsHerzog}, Section 3).
\end{proof}

\begin{lemma}
  Let $X$ be a Gorenstein scheme, which is projective over
  $\boldsymbol{k}$. If $E$ is a coherent sheaf which has a finite resolution
  by injective quasi-coherent sheaves, then $E$ has a finite resolution by
  locally free sheaves. 
\end{lemma}

\begin{proof}
  A finitely generated module over a local Gorenstein ring has finite
  injective dimension if and only if it has finite projective dimension
  \cite{LevinVasco}, Thm. 2.2. Hence, our assumptions imply: $E_{x}$ has
  finite projective dimension for any $x\in X$. Using the Auslander-Buchsbaum
  formula we can bound these projective dimensions by $\dim(X)$. This means
  $$\Ext^{i}(E_{x},-)=0\quad\text{if } i>\dim(X)\;\text{ and }x\in X.$$
  On the other hand, because $X$ is projective, any coherent sheaf has a (not
  necessarily finite) resolution $\cdots\rightarrow L^{-n}\rightarrow L^{-n+1}
  \rightarrow \cdots \rightarrow L^{0}$ by locally free sheaves $L^{i}$ of
  finite rank. If this complex has length greater than $d=\dim(X)$, we replace
  $L^{-d}$ by the kernel $K\subset L^{-d+1}$ of the following map. This
  produces a bounded acyclic complex of coherent sheaves $0\rightarrow
  K\rightarrow L^{-d+1}\rightarrow L^{-d+2} \rightarrow \cdots
  \rightarrow L^{0}\rightarrow E \rightarrow 0$. By standard arguments from
  homological algebra, we obtain
  $$\mathcal{E}xt^{1}(K,-) \cong \mathcal{E}xt^{d+1}(E,-)=0.$$ 
  This implies that all the stalks $K_{x}$ are projective, hence,
  free. Therefore, we have got a finite resolution of $E$ consisting of
  coherent locally free sheaves.
\end{proof}

\begin{corollary}
  A coherent sheaf $E$ on $\boldsymbol{E}$ is spherical, if and only if
  \begin{itemize}
  \item[(i)] $E$ has a finite resolution by locally free sheaves and
  \item[(ii)] $\Hom(E,E)\cong \boldsymbol{k}$. 
  \end{itemize}
\end{corollary}

\begin{example}
  \begin{itemize}
  \item[(a)] The structure sheaf $\boldsymbol{k}(x)$ of a regular point  $x\in
  \boldsymbol{E}$  is spherical.
  \item[(b)] If $L\in\Pic(\boldsymbol{E})$ is a locally free sheaf of rank one,
    we have $\Hom(L,L) \cong H^{0}(L^{\vee}\otimes L) \cong
    H^{0}(\mathcal{O}_{\boldsymbol{E}}) \cong \boldsymbol{k}$. Hence, $L$ is
    spherical. 
  \item[(c)] More generally, any simple vector bundle is spherical.
  \end{itemize}
\end{example}

\begin{remark}
  The structure sheaf $\boldsymbol{k}(y)$ of a singular point $y\in
  \boldsymbol{E}$ does not have a finite locally free resolution, hence it is
  not spherical. 
\end{remark}

Seidel and Thomas carefully define a \emph{twist functor}
$$T_{E}:D^{b}(\boldsymbol{E}) \rightarrow D^{b}(\boldsymbol{E})$$ which is an
exact equivalence, if $E$ is spherical. Basically, for any 
$F\in D^{b}(\boldsymbol{E})$, the object $T_{E}(F)\in D^{b}(\boldsymbol{E})$ 
is the cone over the ``evaluation map''
$$\boldsymbol{R}\Hom(E,F) \otimes E \rightarrow F.$$
More precisely, they replace $E$ and $F$ by injective resolutions $I_{E}$ and
$I_{F}$ which consist of quasi-coherent sheaves and define
$T_{E}(F)$ to be the cone over the map of complexes 
$$\shom(I_{E},I_{F}) \otimes I_{E} \rightarrow I_{F}$$
whose non-zero components are the usual evaluation maps 
$$\Hom(I_{E}^{i},I_{F}^{k+i}) \otimes I_{E}^{i} \rightarrow I_{F}^{k+i}.$$
Here, $\shom(I_{E},I_{F})$ denotes the complex of vector spaces which has
$\oplus_{i}\Hom(I_{E}^{i},I_{F}^{k+i})$ at place $k$. If $E, F$ are coherent
sheaves, the $k$-th cohomology of $\shom(I_{E},I_{F})$ is $\Ext^{k}(E,F)$. 

Because the objects in $I_{F}$ are injective, basic homological algebra shows
that the (genuine) evaluation map
$$\shom(E,I_{F}) \otimes E \rightarrow I_{F}$$
is quasi-isomorphic to the above map 
(see \cite{SeidelThomas}, Prop. 2.6 and Lemma 3.2).  
Now, $\boldsymbol{R}\Hom(E,F)$ is the complex $\shom(E,I_{F})$, hence this map
defines the above ``evaluation map'' as a morphism in $D^{b}(\boldsymbol{E})$.

Most of our computations will be based on the following lemma.

\begin{lemma}\label{computeT}
  Suppose $E, F$ are coherent sheaves on $\boldsymbol{E}$, $E$ is locally free
  and $\Ext^{1}(E,F)=0$. Then:
  $$T_{E}(F) \cong \cone(\Hom(E,F)\otimes E \xrightarrow{\mathrm{ev}} F).$$
\end{lemma}

\begin{proof}
  According to the definition of the twist functor $T_{E}$, the object
  $T_{E}(F)$ is quasi-isomorphic to the cone over the evaluation map 
  $\shom(E,I_{F}) \otimes E \rightarrow I_{F}$, where $F\rightarrow I_{F}$ is
  an injective resolution. Because $\boldsymbol{E}$ is a curve, our assumption
  implies that the natural inclusion $\Hom(E,F) \rightarrow \shom(E,I_{F})$ is
  a quasi-isomorphism. Because $E$ is locally free, this is true for the
  inclusion $\Hom(E,F)\otimes E  \rightarrow \shom(E,I_{F})\otimes E$ as well.
  The claim follows now from the commutative diagram:
  $$\begin{CD}
    \Hom(E,F)\otimes E       @>{\mathrm{ev}}>> F\\
    @VVV                                       @VVV\\
    \shom(E,I_{F})\otimes E  @>{\mathrm{ev}}>> I_{F},
  \end{CD}$$
  in which the horizontal maps are ordinary evaluation maps and the vertical
  maps are quasi-isomorphisms.
\end{proof}

\begin{proposition}
  If $E\in D^{b}(\boldsymbol{E})$ is a spherical object which is isomorphic to
  a bounded complex of locally free sheaves, the twist functor $T_{E}$ is
  isomorphic to the Fourier-Mukai functor $\Phi_{\mathcal{P}}$, whose kernel
  $\mathcal{P}$ is the cone of the natural map
  $\boldsymbol{R}\mathcal{H}om(\pi_{2}^{\ast}E, 
  \pi_{1}^{\ast}E) \rightarrow \mathcal{O}_{\Delta}$. 
\end{proposition}

The proof of \cite{SeidelThomas}, Lemma 3.2, carries over to our situation
literally. 

\begin{example}
  If $E=\mathcal{O}_{\boldsymbol{E}}$ we obtain
  $\boldsymbol{R}\mathcal{H}om(\pi_{2}^{\ast}E, \pi_{1}^{\ast}E) \cong
  \mathcal{O}_{\boldsymbol{E} \times \boldsymbol{E}}$ and $\mathcal{P} \cong
  \mathcal{I}_{\Delta}[1]$, where $\mathcal{I}_{\Delta} \subset
  \mathcal{O}_{\boldsymbol{E} \times \boldsymbol{E}}$ is the ideal sheaf of the
  diagonal $\Delta\subset \boldsymbol{E} \times \boldsymbol{E}$. This implies:
  $$T_{\mathcal{O}} \cong \Phi_{\mathcal{I}_{\Delta}[1]}.$$
\end{example}

\begin{example}
  If $\boldsymbol{k}(x)$ is the structure sheaf of a regular point
  $x\in\boldsymbol{E}$ and $\mathcal{O}_{\boldsymbol{E}}(x)$ is the locally
  free sheaf which corresponds to the Cartier divisor $x$ on $\boldsymbol{E}$,
  \cite{SeidelThomas}, 3.11, shows: $$T_{\boldsymbol{k}(x)}(F) \cong F\otimes
  \mathcal{O}_{\boldsymbol{E}}(x)$$ for any $F\in D^{b}(\boldsymbol{E})$. In
  particular, $T_{\boldsymbol{k}(x)}$ is isomorphic to the Fourier-Mukai
  transform whose kernel is the sheaf
  $\pi_{2}^{\ast}\mathcal{O}_{\boldsymbol{E}}(x)\otimes\mathcal{O}_{\Delta}$.
\end{example}

We shall use the following lemma in the proof of theorem \ref{main} below.  

\begin{lemma}\label{computation}
  Let $x,y\in \boldsymbol{E}$ be closed points and suppose $x$ is a regular
  point. By $\mathcal{I}_{y}\subset\mathcal{O}$ we denote the ideal sheaf of
  $y$. Then there are isomorphisms: 
  \begin{eqnarray}
  T_{\mathcal{O}}(\boldsymbol{k}(y))&\cong &\mathcal{I}_{y}[1]  \notag\\
  T_{\mathcal{O}}(\mathcal{O}(x))&\cong &\boldsymbol{k}(x)  \notag\\
  T_{\mathcal{O}}(\mathcal{O})&\cong &\mathcal{O}.  \notag
  \end{eqnarray}
  If $\boldsymbol{E}$ is singular with singular point $s$, normalisation
  $n:\mathbb{P}^{1}\rightarrow \boldsymbol{E}$ and $\widetilde{\mathcal{O}} =
  n_{\ast}(\mathcal{O}_{\mathbb{P}^{1}})$, we have:
  $$T_{\mathcal{O}}(\widetilde{\mathcal{O}}) \cong \boldsymbol{k}(s).$$
\end{lemma}

\begin{proof}
  The first two isomorphisms are obtained from lemma \ref{computeT} with
  $E=\mathcal{O}$, $F=\boldsymbol{k}(y)$ and $F=\mathcal{O}(x)$ respectively.

  To compute $T_{\mathcal{O}}(\mathcal{O})$, we use an injective resolution
  $\mathcal{O} \rightarrow I_{\mathcal{O}}$ of $\mathcal{O}$ and the exact
  sequence of complexes which is obtained from the definition of the mapping
  cone: 
  $$0 \rightarrow I_{\mathcal{O}} \rightarrow 
      T_{\mathcal{O}}(\mathcal{O}) \rightarrow  
      \shom(\mathcal{O},I_{\mathcal{O}}) \otimes \mathcal{O}[1] \rightarrow  
    0$$
  and its exact cohomology sequence
  $$\begin{CD}
    0 @>>> H^{-1}(T_{\mathcal{O}}(\mathcal{O}))
      @>>> \Hom(\mathcal{O},\mathcal{O})\otimes\mathcal{O}
      @>{\delta}>> \mathcal{O} @>>>\\
      @>>> H^{0}(T_{\mathcal{O}}(\mathcal{O}))
      @>>> \Ext^{1}(\mathcal{O},\mathcal{O})\otimes\mathcal{O}
      @>>> 0.
  \end{CD}$$  
  The connecting homomorphism $\delta$ is in fact the $H^{0}$ of the
  evaluation map $\shom(\mathcal{O},I_{\mathcal{O}}) \otimes \mathcal{O}
  \rightarrow I_{\mathcal{O}}$ whose cone is the complex
  $T_{\mathcal{O}}(\mathcal{O})$. Hence, $\delta$ is the evaluation map
  $\Hom(\mathcal{O},\mathcal{O})\otimes\mathcal{O} \rightarrow \mathcal{O}$
  which is an isomorphism. This implies: $T_{\mathcal{O}}(\mathcal{O})$ is
  isomorphic to a complex which is concentrated in degree zero and
  $H^{0}(T_{\mathcal{O}}(\mathcal{O})) \cong
  \Ext^{1}(\mathcal{O},\mathcal{O})\otimes\mathcal{O} \cong \mathcal{O}$,
  which is the claim.

  Finally, if $\boldsymbol{E}$ is singular, we have
  $H^{0}(\widetilde{\mathcal{O}}) \cong H^{0}(\mathcal{O}_{\mathbb{P}^{1}})
  \cong \boldsymbol{k}$ and $H^{j}(\widetilde{\mathcal{O}}) = 0$ for $j\ne
  0$. The evaluation map $H^{0}(\widetilde{\mathcal{O}}) \otimes \mathcal{O}
  \rightarrow \widetilde{\mathcal{O}}$ is the canonical map $\mathcal{O}
  \rightarrow \widetilde{\mathcal{O}}$. It is injective with cokernel
  $\boldsymbol{k}(s)$. Hence, by lemma \ref{computeT}, we obtain 
  $T_{\mathcal{O}}(\widetilde{\mathcal{O}})\cong\boldsymbol{k}(s)$.
\end{proof}

The key ingredient for the proof of the theorem below is the following result,
versions of which can be found in \cite{BridgeMaciocia}, Section 3.3, and
\cite{FMW}, Lemma 0.3. 

\begin{lemma}\label{identity}
  If $X$ is a projective variety and $\mathbb{G}:D^{b}(X)
  \rightarrow D^{b}(X)$ an integral transform, then the following conditions
  are equivalent: 
  \begin{itemize}
  \item[(i)] $\mathbb{G}$ is isomorphic to the identity functor;
  \item[(ii)] $\mathbb{G}(\mathcal{O}_{X}) \cong
  \mathcal{O}_{X} \;\; \text{ and for all }\; 
  x\in X:\;\; \mathbb{G}(\boldsymbol{k}(x)) \cong \boldsymbol{k}(x).$
  \end{itemize}
\end{lemma}

\begin{proof}
  Denote by $P\in D^{b}(X\times X)$ the kernel of
  $\mathbb{G}$, i.e.\/ $\mathbb{G}(E) \cong \boldsymbol{R}\pi_{2\ast}(P\dtens
  \pi_{1}^{\ast} E)$ for any $E\in D^{b}(X)$. With
  $E=\boldsymbol{k}(x)$ we obtain $\boldsymbol{R}\pi_{2\ast}(P\dtens
  \pi_{1}^{\ast} \boldsymbol{k}(x)) \cong \boldsymbol{L}i_{x}^{\ast}P$, where
  $i_{x}: X\rightarrow X \times X$
  denotes the embedding which satisfies $\pi_{2}\circ i_{x} =
  \Id_{X}$ and $\pi_{1}\circ i_{x}$ is the constant map with image
  $x$. Our assumption implies $\boldsymbol{L}i_{x}^{\ast}P \cong
  \boldsymbol{k}(x)$  is a sheaf for any $x\in X$. By \cite{BridgelandEquiv},
  Lemma 4.3, this implies that $P$ is a sheaf which is $\pi_{1}$-flat. In
  particular, $\boldsymbol{L}i_{x}^{\ast}P \cong 
  i_{x}^{\ast}P \cong \boldsymbol{k}(x)$.

  Because $X$ is projective, we can choose a very ample line
  bundle $A\in\Pic(X)$. Hence, for large positive $m$, the
  canonical mapping
  \begin{equation}
    \label{surjection}\pi_{1}^{\ast}\pi_{1\ast}(P\otimes
    \pi_{2}^{\ast}A^{\otimes m})           \rightarrow
    P\otimes \pi_{2}^{\ast}A^{\otimes m}
  \end{equation}
  is surjective \cite{Hartshorne}, Theorem III.8.8. Observe,
  $h^{j}(i_{x}^{\ast}(P\otimes \pi_{2}^{\ast}A^{\otimes m})) =
  h^{j}(i_{x}^{\ast}(P)\otimes A^{\otimes m}) = h^{j}(\boldsymbol{k}(x)\otimes
  A^{\otimes m}) = h^{j}(\boldsymbol{k}(x)) =1$, if $j=0$ and zero, if
  $j>0$. This implies: $L_{m}:= \pi_{1\ast}(P\otimes \pi_{2}^{\ast}A^{\otimes
    m})$ is locally free of rank one \cite{Hartshorne}, Corollary
  III.12.9. From (\ref{surjection}) we obtain a surjection 
  $$\mathcal{O}_{X\times X} \rightarrow P\otimes
  \pi_{1}^{\ast}L_{m}^{\vee} \otimes \pi{_{2}^{\ast} A^{\otimes m}}.$$
  Because $i_{x}^{\ast}(P\otimes \pi_{1}^{\ast}L_{m}^{\vee} \otimes
  \pi_{2}^{\ast} A^{\otimes m}) \cong \boldsymbol{k}(x)$, there exists a
  unique morphism  $\varphi:X \rightarrow \Hilb^{1}(X)$ which
  satisfies $\varphi^{\ast}(\mathcal{U}) \cong P\otimes
  \pi_{1}^{\ast}L_{m}^{\vee} \otimes \pi_{2}^{\ast} A^{\otimes m}$, where
  $\mathcal{O}_{X \times \Hilb^{1}(X)} \rightarrow \mathcal{U}$ denotes the
  universal quotient sheaf on $\Hilb^{1}(X)$. More precisely,
  $(\Id \times \varphi)^{\ast}(\mathcal{O}_{X \times
    \Hilb^{1}(X)} \rightarrow \mathcal{U})$ is isomorphic to the
  above surjection. In our situation, $\Hilb^{1}(X)
  \cong X$ and $\mathcal{U}\cong \mathcal{O}_{\Delta}$ and the
  morphism must be $\varphi=\Id_{X}$. Hence, $P\otimes
  \pi_{1}^{\ast}L_{m}^{\vee} \otimes \pi_{2}^{\ast} A^{\otimes m} \cong
  \mathcal{O}_{\Delta}$, that is $P \cong \mathcal{O}_{\Delta} \otimes
  \pi_{1}^{\ast}(L_{m}\otimes A^{\otimes -m})$. From $\mathbb{G}(\mathcal{O})
  \cong \mathcal{O}$ we obtain $P\cong \mathcal{O}_{\Delta}$ which means
  $\mathbb{G} \cong \Id$, as desired.
\end{proof}

In the sequel we fix a regular point\label{pnull} $p_{0}\in\boldsymbol{E}$. 
Because $\boldsymbol{E}$ is irreducible, there is the structure of an abelian
group on its regular part and we can choose this structure in such a way that
$p_{0}$ will be the neutral element in this group structure. 
More specifically, the choice of $p_{0}$ allows to define a bijection between
the smooth points on $\boldsymbol{E}$ and the Picard group of invertible
sheaves of degree zero on $\boldsymbol{E}$ which sends a point
$x\in\boldsymbol{E}_{\text{reg}}$ to the sheaf $\mathcal{O}(x-p_{0})$. 
Via this isomorphism, the involution $\mathcal{L}\mapsto \mathcal{L}^{\vee}$ 
on the Picard group corresponds to an involution $i$ on the regular part of
$\boldsymbol{E}$. 
It extends to an involution $i:\boldsymbol{E}\rightarrow\boldsymbol{E}$ 
which fixes the singular point. 
If $x\in\boldsymbol{E}$ is a regular point, $i(x)$ is characterised by the
linear equivalence $i(x)=2p_{0}-x$. 

\begin{definition}
  We define $\mathbb{F}:= T_{\boldsymbol{k}(p_{0})}
  T_{\mathcal{O}_{\boldsymbol{E}}} T_{\boldsymbol{k}(p_{0})}:
  D^{b}(\boldsymbol{E})\rightarrow D^{b}(\boldsymbol{E})$.
\end{definition}

\begin{remark}\label{braid}
  Since $\mathcal{O}_{\boldsymbol{E}}$ and $\boldsymbol{k}(p_{0})$ are
  spherical objects, the functor $\mathbb{F}$ is an exact equivalence.
  Because the two spherical objects $\mathcal{O}_{\boldsymbol{E}}$ and
  $\boldsymbol{k}(p_{0})$ satisfy
  $\Hom^{\ast}(\mathcal{O}_{\boldsymbol{E}},\boldsymbol{k}(p_{0})) \cong
  \boldsymbol{k}$, they form an $A_{2}$ configuration in the sense of
  \cite{SeidelThomas}. From \cite{SeidelThomas}, Prop. 2.13, we obtain
  $$T_{\boldsymbol{k}(p_{0})} T_{\mathcal{O}_{\boldsymbol{E}}}
  T_{\boldsymbol{k}(p_{0})} \cong
  T_{\mathcal{O}_{\boldsymbol{E}}} T_{\boldsymbol{k}(p_{0})} 
  T_{\mathcal{O}_{\boldsymbol{E}}}.$$
\end{remark}

\begin{remark}
  Using the description of $T_{\mathcal{O}}$ as a Fourier-Mukai transform, we
  obtain $\mathbb{F}(E) \cong
  \boldsymbol{R}\pi_{2\ast}(\mathcal{I}_{\Delta}[1] \dtens \pi_{1}^{\ast}(E
  \dtens \mathcal{O}(p_{0}))) \dtens \mathcal{O}(p_{0}) \cong
  \boldsymbol{R}\pi_{2\ast}((\mathcal{I}_{\Delta}[1] \otimes
  \pi_{1}^{\ast}\mathcal{O}(p_{0}) \otimes
  \pi_{2}^{\ast}\mathcal{O}(p_{0})) \dtens \pi_{1}^{\ast}E)$. This is the
  Fourier-Mukai transform with kernel $\mathcal{I}_{\Delta} \otimes
  \pi_{1}^{\ast}\mathcal{O}(p_{0}) \otimes
  \pi_{2}^{\ast}\mathcal{O}(p_{0})[1]$. If $\boldsymbol{E}$ is smooth, up to
  the shift, this is the dual of a Poincar\'e line bundle. From
  \cite{BondalOrlovSemi}, Lemma 1.2, we obtain that, in the smooth case, our
  functor $\mathbb{F}$ is a quasi-inverse of the usual Fourier-Mukai functor
  as studied by Mukai \cite{Mukai}. This explains why the sign of the shift in
  theorem \ref{main} is different from Mukai's.
\end{remark}

\begin{theorem}\label{main}
  It holds $\mathbb{F}\circ \mathbb{F} \cong i^{\ast}[1]$.
  Consequently, $\mathbb{F}^{4} \cong [2]$.
\end{theorem}

\begin{proof}
  Because $\mathbb{F}\circ \mathbb{F}$ and $i^{\ast}[1]$ are Fourier-Mukai
  transforms, by Lemma \ref{identity} we just need to show that these two
  functors give isomorphic objects, if applied to the structure sheaf of
  $\boldsymbol{E}$ and to structure sheaves  of closed  points on
  $\boldsymbol{E}$. 

  Let us compute $\mathbb{F}(\mathbb{F}(\mathcal{O}_{\boldsymbol{E}}))$. From
  $T_{\mathcal{O}}(\mathcal{O})\cong \mathcal{O}$ and remark \ref{braid} we
  see: $\mathbb{F}(\mathbb{F}(\mathcal{O}_{\boldsymbol{E}})) \cong
  T_{\boldsymbol{k}(p_{0})} T_{\mathcal{O}} T_{\mathcal{O}}
  T_{\boldsymbol{k}(p_{0})} (\mathcal{O}_{\boldsymbol{E}})$.
  Hence, and because $i^{\ast}\mathcal{O} \cong \mathcal{O}$, we just need to
  show $T_{\mathcal{O}}T_{\mathcal{O}}
  (\mathcal{O}_{\boldsymbol{E}}(p_{0})) \cong
  \mathcal{O}_{\boldsymbol{E}}(-p_{0})[1]$, but this follows from lemma
  \ref{computation}.

  Next, we compute $\mathbb{F}(\mathbb{F}(\boldsymbol{k}(x)))$,
  where $x\in\boldsymbol{E}$ is a closed point. Again, from
  $T_{\boldsymbol{k}(p_{0})}(\boldsymbol{k}(x)) \cong \boldsymbol{k}(x)$ and
  remark \ref{braid} we deduce: $\mathbb{F}(\mathbb{F}(\boldsymbol{k}(x))
  \cong T_{\mathcal{O}} T_{\boldsymbol{k}(p_{0})} T_{\boldsymbol{k}(p_{0})}
  T_{\mathcal{O}}(\boldsymbol{k}(x))$. If $x$ is a regular point, lemma
  \ref{computation} shows that this is isomorphic to
  $\boldsymbol{k}(i(x))[1]$. Here we used $\mathcal{O}(2p_{0}-x) \cong
  \mathcal{O}(i(x))$, which follows from the definition of $i$. 

  Finally, if $x=s\in\boldsymbol{E}$ is the singular point, we have:
  \begin{eqnarray}  
   T_{\mathcal{O}} T_{\boldsymbol{k}(p_{0})} T_{\boldsymbol{k}(p_{0})}
      T_{\mathcal{O}}(\boldsymbol{k}(s))          &\cong&\notag\\
   T_{\mathcal{O}} T_{\boldsymbol{k}(p_{0})} T_{\boldsymbol{k}(p_{0})}
      (\widetilde{\mathcal{O}(-2)}[1])           &\cong&\notag\\
   T_{\mathcal{O}}  (\widetilde{\mathcal{O}}[1]) &\cong&\notag
   \boldsymbol{k}(s)[1].
  \end{eqnarray}
  Here we used $\mathcal{I}_{s}\cong\widetilde{\mathcal{O}(-2)}\cong
  n_{\ast}\mathcal{O}(-2)$ and $n_{\ast}\mathcal{O}(-2) \otimes
  \mathcal{O}(2p_{0}) \cong n_{\ast}(\mathcal{O}(-2)\otimes
  n^{\ast}\mathcal{O}(2p_{0})) \cong n_{\ast}\mathcal{O}\cong
  \widetilde{\mathcal{O}}$. Because $i(s)=s$ and $i^{\ast}\boldsymbol{k}(s)
  \cong \boldsymbol{k}(s)$ the proof is complete.
\end{proof}

\begin{remark}
  As in the smooth case, we can use this result to obtain an action of the
  group $\widetilde{\SL}(2,\mathbb{Z})$ on the derived category
  $D^{b}(\boldsymbol{E})$. A presentation of $\widetilde{\SL}(2,\mathbb{Z})$
  is given by (see \cite{SeidelThomas}, section 3d):
  $$\langle A,B,T \mid ABA=BAB,\; (AB)^{6}=T^{2},\; [A,T]=[B,T]=1 \rangle.$$
  This is a central extension of $\SL(2,\mathbb{Z})$ by $\mathbb{Z}$, where
  the normal subgroup is generated by $T$ and the projection to
  $\SL(2,\mathbb{Z})$ sends the generators $A,B$ and $T$ to 
  $$
  \begin{pmatrix}
    1&1\\0&1
  \end{pmatrix},\quad
  \begin{pmatrix}
    1&0\\-1&1
  \end{pmatrix} \;\text{and}\quad
  \begin{pmatrix}
    1&0\\0&1
  \end{pmatrix}\;\text{respectively.}
  $$
  The action of $\widetilde{\SL}(2,\mathbb{Z})$ on $D^{b}(\boldsymbol{E})$ is
  obtained by letting the generators $A,B$ and $T$ act as $T_{\mathcal{O}}$,
  $T_{\boldsymbol{k}(p_{0})}$ and as the translation functor $[1]$. 
\end{remark}

For later use, we formulate a consequence of lemma \ref{computeT} explicitly.

\begin{remark}\label{computeF}
  Suppose $E$ is a coherent sheaf on $\boldsymbol{E}$ which satisfies
  $H^{1}(E(p_{0}))=0$. If the evaluation map $$\text{ev}:H^{0}(E(p_{0}))
  \otimes \mathcal{O} \rightarrow E(p_{0})$$
  \begin{itemize}
  \item[(i)] is injective, then $\mathbb{F}(E) \cong \coker(\text{ev})\otimes
    \mathcal{O}(p_{0})$;
  \item[(ii)] is surjective, then $\mathbb{F}(E) \cong \ker(\text{ev})[1]
    \otimes \mathcal{O}(p_{0})$.
  \end{itemize}
\end{remark}

Recall that the degree of a torsion free sheaf $E$ on $\boldsymbol{E}$ is by
definition $\deg(E)=\chi(E)=h^{0}(E)-h^{1}(E)$. Such a sheaf $E$ is called
semi-stable if for any subsheaf $F\subset E$ with $0< \rk(F) < \rk(E)$ it
holds: $$\frac{\deg(F)}{\rk(F)} \le \frac{\deg(E)}{\rk(E)}.$$
Therefore, any torsion free sheaf of rank one is automatically semi-stable. 
If $E$ is semi-stable torsion free and $d>0$ then:
\begin{align*}
  \deg(E)&=d, &\Rightarrow\quad h^{0}(E)&=d, &h^{1}(E)&=0;\\
  \deg(E)&=-d, &\Rightarrow\quad h^{0}(E)&=0, &h^{1}(E)&=d.
\end{align*}

With the aid of theorem \ref{main}, we are able to give a new proof of a
theorem of T.\/ Teodorescu \cite{Teodorescu}. An alternate proof can also be
found in \cite{FMmin}, Cor.~1.2.9.

\begin{theorem}\label{ss}
  For any semi-stable torsion free sheaf $E$ on $\boldsymbol{E}$ of degree
  zero, the evaluation map
  $$\text{ev}:H^{0}(E(p_{0})) \otimes \mathcal{O} \rightarrow E(p_{0})$$ 
  is a monomorphism with cokernel of rank zero.\\
  The functor which sends the sheaf $E$ to the cokernel of this evaluation
  map is the restriction of the functor $\mathbb{F}$. It is an exact
  equivalence between the category of semi-stable torsion free sheaves of
  degree zero and the category of coherent torsion sheaves on $\boldsymbol{E}$.
\end{theorem}

\begin{proof}
  The injectivity of the evaluation map was shown in \cite{FMW}, Theorem 1.2.
  Because $E$ is semi-stable, $E(p_{0})$ is semi-stable as well. From 
  $\deg(E)=0$ we obtain $\deg(E(p_{0}))=\rk(E)$. This implies
  $H^{1}(E(p_{0}))=0$. Hence, by remark \ref{computeF}, $\mathbb{F}(E) \cong
  \coker(\text{ev}) \otimes \mathcal{O}(p_{0})$. Since $h^{0}(E(p_{0})) =
  \deg(E(p_{0})) = \rk(E)$, the cokernel of the evaluation map is a torsion
  sheaf. This implies $\coker(\text{ev}) \cong \coker(\text{ev}) \otimes
  \mathcal{O}(p_{0}) \cong \mathbb{F}(E)$.
  
  From theorem \ref{main} we know that $\mathbb{F}$ is an exact
  equivalence with quasi-inverse $i^{\ast}\circ \mathbb{F}[-1]$. In order to
  complete the proof we just need to show for any torsion sheaf $T$ on
  $\boldsymbol{E}$ that $\mathbb{F}(T)[-1]$ is a semi-stable torsion free sheaf
  of degree zero. By remark \ref{computeF} again, we obtain
  $\mathbb{F}(T)[-1] \cong K(p_{0})$, where $K$ is the kernel of the evaluation
  map, sitting in the exact sequence 
  $$0 \rightarrow K \rightarrow H^{0}(T)\otimes \mathcal{O} \rightarrow T
  \rightarrow 0.$$ This implies that $K$ is torsion free and of rank
  $d=h^{0}(T)$ and degree $-d$. Hence, $\deg(K(p_{0})) = 0$.

  In order to prove
  semi-stability of K, we proceed by induction on $d$. If $d=1$ the rank of
  $K(p_{0})$ is one and $K(p_{0})$ is semi-stable. If $d>1$, we suppose
  $\mathbb{F}(T')[-1]$ is semi-stable for any torsion sheaf $T'$ with
  $h^{0}(T')<d$. There exists a torsion subsheaf $T'\subset T$ with
  $0<h^{0}(T')<h^{0}(T)=d$. This gives an exact commutative diagram
  $$\begin{CD}
    0 @>>> H^{0}(T')\otimes \mathcal{O} @>>> H^{0}(T)\otimes \mathcal{O} @>>>
    H^{0}(T'')\otimes \mathcal{O} @>>> 0\\
      &&    @VVV    @VVV @VVV\\
    0 @>>> T' @>>> T @>>> T'' @>>> 0
  \end{CD}$$
  where $T''$ is another torsion sheaf with $h^{0}(T'')<d$. The vertical maps
  are evaluation maps and their kernels are the torsion free sheaves whose
  twist with $\mathcal{O}(p_{0})$ gives the value of $\mathbb{F}[-1]$ on the
  corresponding torsion sheaf. The snake lemma provides therefore an exact
  sequence of torsion free sheaves
  $$
  \begin{CD}
    0 @>>> \mathbb{F}(T')[-1] @>>> \mathbb{F}(T)[-1] @>>> \mathbb{F}(T'')[-1]
    @>>> 0 
  \end{CD}$$ in which two members are semi-stable by our assumption. Since all
  three sheaves are of degree zero, the middle one is semi-stable as well. 
\end{proof}

Because an indecomposable torsion sheaf on $\boldsymbol{E}$ is supported at
precisely one point, any indecomposable semi-stable torsion free sheaf $E$ of
degree zero has the property that $\supp(\mathbb{F}(E))$ is one point. If
$\supp(\mathbb{F}(E)) = \{x\}$, Friedman, Morgan and Witten \cite{FMW} call
the sheaf $E$ to be \emph{concentrated} at $x$. 

In section \ref{sec:semistable} we shall use theorem \ref{ss} as our main tool
to give an explicit description of \emph{all} semi-stable torsion free sheaves
of degree zero (of finite rank) on a nodal Weierstra{\ss} cubic
$\boldsymbol{E}$.
This extends results of R.\/ Friedman, J.\/ Morgan \cite{FM} who described all
semi-stable locally free sheaves of degree zero on $\boldsymbol{E}$. With our
methods we obtain a new proof of their result.

The following two sections provide the necessary background knowledge about
torsion free sheaves on $\boldsymbol{E}$ as well as torsion sheaves on
$\boldsymbol{E}$ which are supported at the singularity $s\in\boldsymbol{E}$,
which is supposed to be a node. 

\section{Torsion free sheaves}
\label{sec:tf}

Throughout this section, $\boldsymbol{E}$ denotes a Weierstra{\ss}
cubic with one node (ordinary double point) $s\in\boldsymbol{E}$. 
The classification of indecomposable vector bundles in the smooth case is now
classical, see \cite{Atiyah}. 
In \cite{DrozdGreuel}, Drozd and Greuel applied methods from representation
theory to the classification of indecomposable torsion free sheaves on
$\boldsymbol{E}$. It is possible to deduce from their results that the sheaves
described below comprise all indecomposable torsion free sheaves on
$\boldsymbol{E}$. They describe these sheaves by gluing fibres of sheaves on
the normalisation, a method which was already used by Seshadri \cite{Seshadri}.
We prefer a description of such sheaves as direct images of
line bundles under finite morphisms. Friedman, Morgan \cite{FM}, Section 2.4,
and Teodorescu \cite{Teodorescu} used similar methods to describe stable and
semi-stable vector bundles.

\subsection{Unipotent vector bundles} 

As in \cite{Atiyah}, for any integer $m\ge 1$ an indecomposable vector bundle
$\mathcal{F}_{m}$ of rank $m$ which satisfies $\mathcal{F}_{m}^{\vee} \cong
\mathcal{F}_{m}$ and $\Ext^{1}(\mathcal{F}_{m}, \mathcal{O}) \cong
\boldsymbol{k}$ is defined as follows: $\mathcal{F}_{1}\cong \mathcal{O}$ and
for $m\ge 1$ the bundle $\mathcal{F}_{m+1}$ is the unique one which appears in
a non-split extension
$$0\rightarrow \mathcal{O} 
   \rightarrow \mathcal{F}_{m+1} 
   \rightarrow \mathcal{F}_{m} 
   \rightarrow 0.$$

\subsection{\'Etale covers}\label{etalecovers} 
Let $\boldsymbol{E_{n}}$ be a cycle of $n$ rational
curves. This means, $\boldsymbol{E_{n}}$ is a connected reduced curve with $n$
ordinary double points whose normalisation has $n$ connected components $D_{1},
D_{2},\ldots,D_{n}$ each of which is isomorphic to $\mathbb{P}^{1}$. 
In addition, the singularities $s_{1},\ldots,s_{n}$ are the points where these
components are glued together. We choose notation in such a way that $s_{\nu}$
is the point where $D_{\nu}$ and $D_{\nu+1}$ are glued together. Here and
below, we use cyclic subscripts, i.e.\/ $D_{n+k}=D_{k}$ for any integer $k$.
Up to isomorphism, there is only one such curve for any $n\ge 1$. It is a
reduced curve of arithmetic genus one.  
If $n=1$ we obtain $\boldsymbol{E_{1}}= \boldsymbol{E}$. 

Once and for all, we fix projective coordinates $(x_{\nu}:y_{\nu})$ on each
$D_{\nu} \cong \mathbb{P}^{1}$ in such a way that $s_{\nu-1}''=(1:0)$ and
$s_{\nu}'=(0:1)$ are the points on $D_{\nu}$ which are mapped to $s_{\nu-1} \in
\boldsymbol{E_{n}}$ and $s_{\nu} \in \boldsymbol{E_{n}}$ respectively.

Furthermore, we fix a normalisation map $\pi_{1}: D_{1}\rightarrow
\boldsymbol{E_{1}}$ which sends the point with coordinates $(1:-1)$ to the
point $p_{0}\in\boldsymbol{E}$ (see p.\/ \pageref{pnull}). The morphism
$\coprod D_{\nu} \rightarrow D_{1}$ which is the identity map on each component
with respect to the coordinates chosen above, descends to an \'etale morphism 
$\pi_{n}: \boldsymbol{E_{n}} \rightarrow \boldsymbol{E}$. 
For any $n\ge 1$, this is a cyclic Galois cover of degree $n$. 
Our choices imply that the point $p_{0\nu}\in\boldsymbol{E_{n}}$ which
corresponds to $(1:-1)\in D_{\nu}$ satisfies $\pi_{n}(p_{0\nu}) = p_{0}$.

We obtain local coordinates $x_{\nu}$ and $y_{\nu+1}$ for the two branches of
$\boldsymbol{E_{n}}$ which intersect at $s_{\nu}$. 
Therefore, the completion $R_{\nu}$ of the local ring of $\boldsymbol{E_{n}}$
at $s_{\nu}$ is isomorphic to 
\[
R_{\nu} \cong \local{x_{\nu}}{y_{\nu+1}}.
\]

We adjust these isomorphisms in such a way that
\begin{equation}\label{isom}
\local{x}{y} \xrightarrow{\sim} \local{x_{\nu}}{y_{\nu+1}}
\end{equation}
which sends $x$ to $x_{\nu}$ and $y$ to $y_{\nu+1}$, is the isomorphism
induced by $\pi_{n}$ between the completion $R$ of the local
ring of $\boldsymbol{E}$ at the singular point $s$ and $R_{\nu}$.

\subsection{Line bundles on $\boldsymbol{E_{n}}$}\label{linebundles}
We denote $\boldsymbol{d}p_{0} := \sum_{i=1}^{n} d_{i}p_{0i}$, which is a
divisor on $\boldsymbol{E_{n}}$ supported in the regular locus. 
It is well known that $$\Pic(\boldsymbol{E_{n}}) \cong
\mathbb{Z}^{n} \times \boldsymbol{k}^{\times}.$$ 
We choose these isomorphisms such that the line bundle
$\mathcal{L}=\mathcal{L}(\boldsymbol{d},\lambda)$, which corresponds to
$(\boldsymbol{d},\lambda) = ((d_{1},\ldots,d_{n}),\lambda) \in \mathbb{Z}^{n}
\times \boldsymbol{k}^{\times}$, satisfies:
$$d_{\nu}= \deg(\mathcal{L}|_{D_{\nu}}).$$
We fix notation by the two requirements:
$$
  \mathcal{L}(\boldsymbol{d},1) \cong \mathcal{O}(\boldsymbol{d}p_{0})
  \quad\text{and}\quad
  \mathcal{L}(\boldsymbol{0},\lambda^{n}) \cong
  \pi_{n}^{\ast}(\mathcal{L}(0,\lambda)).
$$
In case $n=1$ all line bundles of degree one are of the form $\mathcal{O}(p)$
with a regular point $p\in\boldsymbol{E}$. Hence, we obtain a bijection
$P:\boldsymbol{k}^{\times} \rightarrow \boldsymbol{E}_{\text{reg}}$ 
which satisfies $P(1)=p_{0}$ and
$\mathcal{L}(1,\lambda) \cong \mathcal{O}(P(\lambda))$. Thus, for any
$(d,\lambda) \in\mathbb{Z} \times  \boldsymbol{k}^{\times}$, one has 
$$\mathcal{L}(d,\lambda) \cong \mathcal{O}(P(\lambda)+(d-1)p_{0}).$$
With this notation, the involution $i:\boldsymbol{E}_{\text{reg}} \rightarrow
\boldsymbol{E}_{\text{reg}}$ is described as $i(P(\lambda)) =
P(\lambda^{-1})$, i.e.\/ $P$ is a homomorphism.

The above conditions fix our choices, because
\begin{align*}
\mathcal{L}(\boldsymbol{d},\lambda) 
&\cong 
\mathcal{L}(\boldsymbol{d},1) \otimes \mathcal{L}(\boldsymbol{0},\lambda) \\
&\cong
\mathcal{O}(\boldsymbol{d}p_{0}) \otimes \pi_{n}^{\ast}(\mathcal{L}(0,\xi))\\
&\cong
\mathcal{O}(\boldsymbol{d}p_{0}) \otimes 
                                  \pi_{n}^{\ast}(\mathcal{O}(P(\xi)-p_{0}))\\
&\cong
\mathcal{O}(\pi_{n}^{-1}(P(\xi)) + (\boldsymbol{d-1})p_{0})
\end{align*}
where $\xi\in\boldsymbol{k}^{\times}$ is arbitrary with $\xi^{n}=\lambda$.

Our choices are made in such a way that $\mathcal{L}(\boldsymbol{d},\lambda)$
can be described by gluing the line bundles $\mathcal{O}(d_{\nu})$. 
Using the coordinates chosen above and the standard local trivialisation of
$\mathcal{O}(d_{\nu})$ on $D_{\nu}$, the line bundle
$\mathcal{L}(\boldsymbol{d},\lambda)$ is obtained by gluing with the identity
at $s_{1},\ldots,s_{n-1}$ but with $\lambda$ at $s_{n}$. This will be made more
precise below. 

Using the standard trivialisations means that the section $x^{a}y^{d-a}$ of the
line bundle $\mathcal{O}(d)$ on $\mathbb{P}^{1}$ has a local description as
$x^{a}$ about the point $(0:1)$ where $x$ is a local coordinate. 
About the point $(1:0)$, where $y$ is a local coordinate, the local description
is $y^{d-a}$. More generally, the section of $\mathcal{O}(d)$, which is given
by a homogeneous polynomial $f(x,y)$ of degree $d$, has local descriptions
$f(x,1)$ about the point $(0:1)$ and $f(1,y)$ about the point $(1:0)$. The
bundles $\mathcal{O}(d_{\nu})$ are glued together to give
$\mathcal{L}(\boldsymbol{d},\lambda)$ in such a way that homogeneous
polynomials $f_{\nu} \in H^{0}(D_{\nu},\mathcal{O}(d_{\nu}))$ represent a
section of $\mathcal{L}(\boldsymbol{d},\lambda)$ if and only if they satisfy
the gluing condition:
\begin{equation}\label{gluing}
  \begin{aligned}
    f_{\nu}(0:1) &= f_{\nu+1}(1:0) &1\le \nu < n\\
    f_{n}(0:1) &= \lambda f_{1}(1:0).
  \end{aligned}
\end{equation}
To see that this matches our conventions, observe that
$\mathcal{L}(\boldsymbol{d},1)$ has a global section which is given by 
$f_{\nu} =(x_{\nu}+y_{\nu})^{d_{\nu}} \in H^{0}(D_{\nu},\mathcal{O}(d_{\nu}))$.
This section vanishes on each component precisely at $p_{0\nu}$, because
$f_{\nu}$ vanishes on $D_{\nu}$ precisely at the point with coordinates
$(1:-1)$. On the other hand, denoting $\boldsymbol{1}:=(1,1,\ldots,1)\in
\mathbb{Z}^{n}$, $\mathcal{L}(\boldsymbol{1},\lambda^{n}) \cong 
\pi_{n}^{\ast}\mathcal{L}(1,\lambda)$ has a global section which is given by
$f_{\nu} = \lambda^{\nu-1}(x_{\nu}+\lambda y_{\nu}) \in
H^{0}(D_{\nu},\mathcal{O}(1))$. This section vanishes on each component
$D_{\nu}$ at the point with coordinates $(\lambda:-1)$. 
This means in particular $P(\lambda)=\pi_{1}(\lambda:-1)$ and $i(P(\lambda)) =
\pi_{1}(-1:\lambda)$, hence $i:\boldsymbol{E}\rightarrow\boldsymbol{E}$ is the
morphism which is induced by the map $(x:y) \mapsto (y:x)$ on the
normalisation $D_{1}\cong \mathbb{P}^{1}$.

\subsection{Vector bundles via \'etale covers} 
$$\mathcal{B}(\boldsymbol{d},m,\lambda) :=
\pi_{n\ast}\mathcal{L}(\boldsymbol{d},\lambda) \otimes \mathcal{F}_{m}$$
is a vector bundle of rank $mn$ and degree $m\sum_{\nu=1}^{n} d_{n}$. 

It is indecomposable if and only if $\boldsymbol{d}$ is non-periodic. This
means that there does not exist an integer $n'<n$ and a vector
$\boldsymbol{e}\in \mathbb{Z}^{n'}$ such that
$\boldsymbol{d}=(\boldsymbol{e},\boldsymbol{e},\ldots,\boldsymbol{e})$.
Equivalently, $\mathcal{L}(\boldsymbol{d},\lambda)$ is not the pull back
$\pi_{n',n}^{\ast}(\mathcal{L(\boldsymbol{e},\lambda')})$ of a line bundle
under an \'etale morphism $\pi_{n',n}:\boldsymbol{E_{n}} \rightarrow
\boldsymbol{E_{n'}}$.

\subsection{Chains of lines} 
Let $\boldsymbol{I_{n}} \subset \boldsymbol{E_{n+1}}$ be the
chain of projective lines which is obtained from $\boldsymbol{E_{n+1}}$ by
removing one component (say $D_{n+1}$). The singularities of
$\boldsymbol{I_{n}}$ are the points $s_{1}, s_{2}, \ldots, s_{n-1}$. We use the
same coordinates and conventions as before, with the exception that the
subscripts are not cyclic. 
We denote the restriction of $\pi_{n+1}$ to $\boldsymbol{I_{n}}$ by 
$p_{n}:\boldsymbol{I_{n}} \rightarrow \boldsymbol{E}$. 
The two smooth points in $p_{n}^{-1}(s)$ are denoted by
$s_{0}$ (corresponding to $(1:0)\in D_{1}$) and $s_{n}$ (corresponding to
$(0:1)\in D_{n}$). If $n=1$ we have
$\boldsymbol{I_{1}}\cong\mathbb{P}^{1}$ and $p_{1}:\boldsymbol{I_{1}}
\rightarrow \boldsymbol{E}$ is the normalisation. In this case, $s_{0}$ and
$s_{1}$ are the two preimages of the singular point $s\in \boldsymbol{E}$. 
We have $\Pic(\boldsymbol{I_{n}}) \cong \mathbb{Z}^{n}$ and
$\mathcal{L}=\mathcal{L}(\boldsymbol{d})\in \Pic(\boldsymbol{I_{n}})$ is
determined by the degrees $d_{\nu}=\deg(\mathcal{L}|_{D_{\nu}})$.

The description of the completed local rings at the singularities $s_{\nu}$
does not differ from the case $\boldsymbol{E_{n+1}}$. If $1\le \nu \le n-1$
the mappings $R \rightarrow R_{\nu}$ between the completed local rings which
are induced by $p_{n}$ have the same description as above. As before,
$R=\local{x}{y}$. 

The completed local rings at the regular points $s_{0}$ and $s_{n}$ are
isomorphic to $R_{0}=\boldsymbol{k}[[y_{1}]]$ and $R_{n} =
\boldsymbol{k}[[x_{n}]]$. We choose these isomorphisms such that the mapping
$R\rightarrow R_{0}$, which is induced by $p_{n}$, sends $x$ to $0$ and $y$ to
$y_{1}$. Similarly, $R \rightarrow R_{n}$ sends $x$ to $x_{n}$ and $y$ to $0$.

\subsection{Non-locally free sheaves} 
$$\mathcal{S}(\boldsymbol{d}):= p_{n\ast}\mathcal{L}(\boldsymbol{d})$$ 
is an indecomposable torsion free sheaf of rank $n$ and degree
$1+\sum_{\nu=1}^{n}d_{\nu}$ on $\boldsymbol{E}$. The sheaf
$\mathcal{S}(\boldsymbol{d})$ is not locally free. 
If $n=1$ and $d\in\mathbb{Z}$, $\mathcal{S}(d)$ is the sheaf which was denoted
$\widetilde{\mathcal{O}(d)}$ in section \ref{sec:fm}.

\subsection{Useful results}

\begin{lemma}\label{cohomology}
  Let $m\ge 1, \lambda\in\boldsymbol{k}^{\times}$ and
  $\boldsymbol{d}\in\mathbb{Z}^{n}$ which satisfies $d_{\nu}\ge 0$ for any 
  $\nu$ and $\sum_{\nu=1}^{n} d_{\nu} > 0$. Then, the following holds:
  \begin{align*}
    h^{0}(\mathcal{B}(\boldsymbol{d},m,\lambda))&= m\sum_{\nu=1}^{n} d_{n} &
    h^{0}(\mathcal{S}(\boldsymbol{d})) &=1+\sum_{\nu=1}^{n} d_{n} \\
    h^{1}(\mathcal{B}(\boldsymbol{d},m,\lambda))&=0 &
    h^{1}(\mathcal{S}(\boldsymbol{d})) &=0
  \end{align*}
\end{lemma}

\begin{proof}
  We defined
  $\mathcal{B}(\boldsymbol{d},1,\lambda) =
  \pi_{n\ast}\mathcal{L}(\boldsymbol{d},\lambda)$. Finiteness of
  $\pi_{n}:\boldsymbol{E_{n}}\rightarrow \boldsymbol{E}$ implies
  $H^{i}(\mathcal{B}(\boldsymbol{d},1,\lambda)) \cong
  H^{i}(\mathcal{L}(\boldsymbol{d},\lambda))$. Let $\eta: \coprod_{\nu=1}^{n}
  D_{\nu} \rightarrow \boldsymbol{E_{n}}$ denote the normalisation. There is an
  exact sequence  
  $$\begin{CD}
    0 @>>> \mathcal{L}(\boldsymbol{d},\lambda)
      @>>> \eta_{\ast}\eta^{\ast}\mathcal{L}(\boldsymbol{d},\lambda)
      @>{\alpha}>> \bigoplus_{\nu=1}^{n} \boldsymbol{k}(s_{\nu})
      @>>> 0
    \end{CD}$$
  Because $H^{0}(\eta_{\ast}\eta^{\ast}\mathcal{L}(\boldsymbol{d},\lambda))
  \cong H^{0}(\eta^{\ast}\mathcal{L}(\boldsymbol{d},\lambda)) \cong
  \bigoplus_{\nu=1}^{n} H^{0}(D_{n},\mathcal{O}(d_{n}))$, 
  $H^{0}(\alpha)=0$ is precisely our gluing condition. With the choices made
  above, an explicit description of $H^{0}(\alpha)$ is the
  following: the $\nu$-th component of $H^{0}(\alpha)(f_{1},\ldots,f_{n})$ is
  \begin{align*}
    f_{\nu}(0:1) &- f_{\nu+1}(1:0)     &\text{if }\;\; \nu<n\\
    f_{n}(0:1)   &- \lambda f_{1}(1:0) &\text{if }\;\; \nu=n.
  \end{align*}
  Because all $d_{n}\geq 0$ and at least one of them is positive, it is now
  easy to see that $H^{0}(\alpha)$ is surjective. Using
  $H^{1}(\eta_{\ast}\eta^{\ast}\mathcal{L}(\boldsymbol{d},\lambda)) \cong
  H^{1}(\eta^{\ast}\mathcal{L}(\boldsymbol{d},\lambda)) \cong
  \bigoplus_{\nu=1}^{n} H^{1}(D_{n},\mathcal{O}(d_{n})) = 0$, the exact
  sequence implies $H^{1}(\mathcal{L}(\boldsymbol{d},\lambda))=0$ and
  $h^{0}(\mathcal{L}(\boldsymbol{d},\lambda)) =
  h^{0}(\oplus\mathcal{O}(d_{\nu})) -n = \sum_{\nu=1}^{n} (d_{\nu}+1) -n =
  \sum_{\nu=1}^{n} d_{\nu}$.
  Using induction and the exact sequences
  $$\begin{CD}
    0 @>>> \mathcal{B}(\boldsymbol{d},1,\lambda)
      @>>> \mathcal{B}(\boldsymbol{d},m+1,\lambda)
      @>>> \mathcal{B}(\boldsymbol{d},m,\lambda)
      @>>> 0
    \end{CD}$$
  the remaining statements about $\mathcal{B}(\boldsymbol{d},m,\lambda)$
  follow. In particular, we have shown
  $\deg(\mathcal{B}(\boldsymbol{d},m,\lambda)) = m\sum_{\nu=1}^{n} 
  d_{n}$ under the assumptions made.

  The sheaves $\mathcal{S}(\boldsymbol{d})$ are studied in the same way. This
  time we use the finite morphism $p_{n}:\boldsymbol{I_{n}}\rightarrow
  \boldsymbol{E}$ and a normalisation $\eta: \coprod_{\nu=1}^{n} D_{\nu}
  \rightarrow \boldsymbol{I_{n}}$. We obtain an exact sequence 
  $$\begin{CD}
    0 @>>> \mathcal{L}(\boldsymbol{d})
      @>>> \eta_{\ast}\eta^{\ast}\mathcal{L}(\boldsymbol{d})
      @>{\beta}>> \bigoplus_{\nu=1}^{n-1} \boldsymbol{k}(s_{\nu})
      @>>> 0
    \end{CD}$$
  The description of $H^{0}(\beta)$ is as above but without the $n$-th
  component involving $\lambda$. The surjectivity is seen easily. We might
  allow all $d_{\nu}$ to be zero here. As above, we conclude with
  $H^{1}(\mathcal{S}(\boldsymbol{d}))=0$ and 
  $h^{0}(\mathcal{S}(\boldsymbol{d})) = h^{0}(\oplus_{\nu=1}^{n}
  \mathcal{O}(d_{\nu})) - (n-1) = 1+\sum_{\nu=1}^{n} d_{\nu}$.
\end{proof}

For $\mathcal{B}(\boldsymbol{d},m,\lambda)$ with arbitrary $\boldsymbol{d}$
this result is contained in \cite{DGK}. We include a proof here, because we
did not show that our definition of the sheaves
$\mathcal{B}(\boldsymbol{d},m,\lambda)$ coincides with theirs.

The following lemma will be used at the end of section \ref{sec:matlis}. 

\begin{lemma}\label{compDual}
  If $\boldsymbol{d}=(d_{1},\ldots,d_{n})\in \mathbb{Z}^{n},
  \lambda\in\boldsymbol{k}^{\times}, m\ge 1$, we have:
  \begin{itemize}
  \item[(i)] $\mathcal{B}(\boldsymbol{d},m,\lambda)^{\vee} \cong
              \mathcal{B}(-\boldsymbol{d},m,\lambda^{-1})$
  \item[(ii)] $\mathcal{S}(\boldsymbol{d})^{\vee} \cong 
               \mathcal{S}(\boldsymbol{\kappa - d})$ with 
              $\boldsymbol{\kappa}=
              \begin{cases}
                (-1,0,\ldots,0,-1) &\text{if } n\ge 2\\
                -2             &\text{if } n= 1.
              \end{cases}$
  \end{itemize}
\end{lemma}

\begin{proof}
  If $f:X\rightarrow \boldsymbol{E}$ is a finite morphism, $F$ a coherent sheaf
  on $X$ and $G$ a coherent sheaf on $\boldsymbol{E}$, there is a natural
  $f_{\ast}\mathcal{O}_{X}$-morphism
  $$f_{\ast}\mathcal{H}om_{X}(F,f^{!}G) \cong 
            \mathcal{H}om_{\boldsymbol{E}}(f_{\ast}F,G).$$
  The coherent $\mathcal{O}_{X}$-module $f^{!}G$ is characterised by the
  isomorphism of $f_{\ast}\mathcal{O}_{X}$-modules $f_{\ast}f^{!}G \cong 
    \mathcal{H}om_{\boldsymbol{E}}(f_{\ast}\mathcal{O}_{X},G).$
    Recall that $f^{!}\omega_{\boldsymbol{E}}$ is a dualising sheaf on $X$ if
    $\omega_{\boldsymbol{E}}$ is one on $\boldsymbol{E}$. In our situation
    $\omega_{\boldsymbol{E}} \cong \mathcal{O}_{\boldsymbol{E}}$ and we obtain
    an isomorphism
    $$f_{\ast}\mathcal{H}om_{X}(F,\omega_{X}) \cong
    \mathcal{H}om_{\boldsymbol{E}}(f_{\ast}F,\mathcal{O}_{\boldsymbol{E}})
    \cong (f_{\ast}F)^{\vee}.$$
    To show (i), we consider $X=\boldsymbol{E_{n}}$ and $f=\pi_{n}$. The claim
    follows now from $\omega_{\boldsymbol{E_{n}}} \cong
    \mathcal{O}_{\boldsymbol{E_{n}}}$, $\mathcal{F}_{m}^{\vee}
    \cong \mathcal{F}_{m}$ and
    $\mathcal{L}(\boldsymbol{d},\lambda)^{\vee} \cong
    \mathcal{L}(-\boldsymbol{d},\lambda^{-1})$ on $\boldsymbol{E_{n}}$.

    For the proof of (ii) we let $X=\boldsymbol{I_{n}}$ and $f=p_{n}$. Now
    $\omega_{\boldsymbol{I_{n}}} \cong \mathcal{L}(\boldsymbol{\kappa})$ and
    the result follows from $\mathcal{L}(\boldsymbol{d})^{\vee} \cong
    \mathcal{L}(-\boldsymbol{d})$ on $\boldsymbol{I_{n}}$. 
\end{proof}

Based on the results collected above, we easily obtain:

\begin{corollary}\label{rankone}
  Any locally free sheaf on $\boldsymbol{E}$ which is isomorphic to
  $\mathcal{L} \otimes \mathcal{F}_{m}$ with
  $\mathcal{L}\in\Pic(\boldsymbol{E})$ is semi-stable.
  For any $d\in\mathbb{Z}$, the rank one torsion free sheaves $\mathcal{S}(d)$
  are semi-stable. 
  On the sheaves of degree zero among those, the functor $\mathbb{F}$
  has the following description: 
  $$\mathbb{F}(\mathcal{S}(-1))\cong \boldsymbol{k}(s) \;\;\text{ and }\;\;
  \mathbb{F}(\mathcal{B}((0),m,\lambda)) \cong
  \mathcal{O}_{P(\lambda)}/\mathfrak{m}_{P(\lambda)}^{m}.$$ 
\end{corollary}

\begin{proof}
  The semi-stability is clear from the definitions. The description of
  $\mathbb{F}$ follows easily from theorem \ref{ss} if the rank is one. If
  $m>1$ we use the extensions which define the vector bundles
  $\mathcal{F}_{m}$ and proceed by induction. 
\end{proof}

In section \ref{sec:semistable} we shall describe the semi-stable torsion free
sheaves of degree zero and arbitrary rank as well as their images under
$\mathbb{F}$. This requires some knowledge about torsion sheaves supported at
the node $s\in\boldsymbol{E}$ which will be collected in section
\ref{sec:tors}.

\section{Torsion Sheaves}
\label{sec:tors}

In this section we study indecomposable modules of finite length over the
complete local ring $R = \local{x}{y}$. Such a module is given by a finite
dimensional vector space $V$ over $\boldsymbol{k}$ and two commuting
$\boldsymbol{k}$-linear nilpotent endomorphisms $X,Y : V\rightarrow V$ which
satisfy $XY=0$. Two modules $(V,X,Y)$ and $(V',X',Y')$ are isomorphic if there
is an isomorphism of $\boldsymbol{k}$-vector spaces $S: V\rightarrow V'$ such
that $SX = X'S, SY = Y'S$. Hence, the classification of $R$-modules of finite
length is equivalent to the classification of pairs of commuting nilpotent 
matrices $(X,Y)$ with $XY=0$ modulo simultaneous conjugation.
This problem was solved by Gelfand and Ponomarev
\cite{GelfandPonomarev} (see also \cite{NazRoi}, \cite{BNRS}). They gave a
complete classification of such modules.

Let us recall their results. There are two types of indecomposable $R$-modules
of finite length: 
the so-called \emph{bands} (forming a family depending on one continuous and
some discrete parameters)  and \emph{strings} (a discrete family).  
In \cite{GelfandPonomarev} strings were called modules of the \emph{first kind}
and bands were called modules of the \emph{second kind}.
For our purposes the following description of indecomposable objects is the
most appropriate. Using the band and string diagrams introduced below, it is
not hard to derive this description from the explicit description which was
given in \cite{BNRS}, see \cite{BurbanDrozd}. 

\subsection*{Bands}
  A \emph{band} $\mathcal{M}(\boldsymbol{q},m,\lambda)$ depends on an
  integer $m \ge 1$, a parameter $\lambda \in \boldsymbol{k}^{\times}$ and a 
  non-periodic sequence of pairs of integers 
  $$\boldsymbol{q} = (n_1, m_1)(n_2, m_2)\dots (n_N, m_N)$$
  with $n_i, m_i \ge 1$ and $N\ge 1$. 
  Its  minimal free resolution  is 
  $$\begin{CD}
    0 @>>> R^{mN} @>{M(\boldsymbol{q},m,\lambda)}>> R^{mN} @>>>
      \mathcal{M}(\boldsymbol{q},m,\lambda) @>>> 0,  
    \end{CD}$$
  where 
  $$M(\boldsymbol{q},m,\lambda) = 
    \begin{pmatrix}
      x^{n_1}I_{m} & y^{m_1}I_{m} & 0          & \dots & 0 \\
        0          & x^{n_2}I_{m} & y^{m_2}I_{m} & \dots & 0 \\
        \vdots     & \vdots     &\ddots      & \ddots& \vdots     \\
        0          & 0          & 0          & \ddots& y^{m_{N-1}}I_{m}\\
        y^{m_N}J_{m}(\lambda) & 0 & 0          & \dots & x^{n_N}I_{m}
     \end{pmatrix}.$$
If $N=1$ this reads as
$$
M((n_{1},m_{1}),m,\lambda) = (x^{n_1}I_{m} + y^{m_1}J_{m}(\lambda)).
$$
Here, by $I_{m}$ we denote the identity matrix of size $m\times m$ and by
$J_{m}(\lambda)$ we denote the Jordan block of size $m\times m$ and with
eigenvalue $\lambda$. 
We obtain the same module, if we place $J_{m}(\lambda)$ as a factor at any
other power of $y$, whereby we ensure that it appears only once in the
matrix. Because $J_{m}(\lambda)^{-1}$ and $J_{m}(\lambda^{-1})$ are similar
matrices, instead we could place $J_{m}(\lambda^{-1})$ at precisely one
position as a factor of a power of $x$.

Any band module has finite injective dimension. A cyclic permutation of the
pairs which constitute $\boldsymbol{q}$ does not change the band module.

\subsection*{Strings}
A \emph{string} module $\mathcal{N}(\boldsymbol{q})$ depends on a
  sequence of integers
$$\boldsymbol{q} = n_0 (m_1, n_1)(m_2, n_2)\dots (m_N, n_N) m_{N+1},$$
where  $n_0, m_{N+1} \ge 0$ and $n_i, m_i \ge 1$ for $ 1\le i \le N$. It is
convenient to write their resolution in the form
$$\begin{CD}
  0 @>>> R^{N+1} @>{N(\boldsymbol{q})}>> \boldsymbol{k}[[y]] 
\oplus R^N \oplus \boldsymbol{k}[[x]]  @>>>  \mathcal{N}(\boldsymbol{q}) @>>>
0, 
\end{CD}$$
where 
$$
N(\boldsymbol{q}) = 
\begin{pmatrix}
y^{n_0}&   0     &  0          &\hdotsfor[2]{1}&   0    \\
x^{m_1}& y^{n_1} &  0          &\hdotsfor[2]{1}&   0    \\
0      & x^{m_2} &  y^{n_2}    &\hdotsfor[2]{1}&\vdots  \\
0      &\hdotsfor[2]{3}                        &   0    \\
\vdots &   0     &  0          &   x^{m_N}     & y^{n_N}\\
0      &   0     &  0          &    0          & x^{m_{N+1}}
 \end{pmatrix}.
$$
If $N=0$ we obtain $\boldsymbol{q}=n()m$ and the above exact sequence becomes 
$$\begin{CD}
  0 @>>> R                          @>{\binom{y^{n}}{x^{m}}}>>
  \boldsymbol{k}[[y]] \oplus \boldsymbol{k}[[x]]  @>>>
  \mathcal{N}(n()m)                               @>>> 0.
\end{CD}$$
If $\boldsymbol{q}=0()0$, this means $\mathcal{N}(0()0) \cong
\boldsymbol{k}(s)$. In all other cases the module has length at least two.
String modules do not have finite injective dimension. 

According to \cite{GelfandPonomarev}, any indecomposable $R$-module of finite
length is isomorphic to one of the bands or strings described above. Hence, any
indecomposable $R$-module of finite length and of finite injective dimension
is isomorphic to a module $\mathcal{M}(\boldsymbol{q},m,\lambda)$. Whereas, if
such a module has infinite injective dimension, it must be isomorphic to a
module $\mathcal{N}(\boldsymbol{q})$. 

\begin{remark}\label{dummy}
  With $\boldsymbol{q}$ as above, the matrix
  $$N'(\boldsymbol{q}) := 
   \begin{pmatrix}
    1      &0      &\hdotsfor[2]{3}                        &   0    \\
    1      &y^{n_0}&   0     &  0          &\hdotsfor[2]{1}&   0    \\
    0      &x^{m_1}& y^{n_1} &  0          &\hdotsfor[2]{1}&   0    \\
    0      &0      & x^{m_2} &  x^{n_2}    &\hdotsfor[2]{1}&\vdots  \\
    0      &0      &\hdotsfor[2]{3}                        &   0    \\
    \vdots &\vdots &   0     &  0          &   x^{m_N}     & y^{n_N}\\
    0      &0      &   0     &  0          &    0          & x^{m_{N+1}}
   \end{pmatrix}
  $$
  defines a linear mapping $R^{N+2}\rightarrow \boldsymbol{k}[[y]] \oplus
  R^{N+1} \oplus \boldsymbol{k}[[x]]$ whose cokernel is isomorphic to
  $\mathcal{N}(\boldsymbol{q})$. To see this, we subtract $y^{n_{0}}$ times the
  first column from column two and add row one to row two in
  $N'(\boldsymbol{q})$. This corresponds to applying $R$-linear automorphisms
  to $R^{N+2}$ and $\boldsymbol{k}[[y]] \oplus R^{N+1} \oplus
  \boldsymbol{k}[[x]]$ and transforms $N'(\boldsymbol{q})$ to:
  $$\widetilde{N}'(\boldsymbol{q}) := 
   \begin{pmatrix}
    0      &y^{n_0}&\hdotsfor[2]{3}                        &   0    \\
    1      &0      &   0     &  0          &\hdotsfor[2]{1}&   0    \\
    0      &x^{m_1}& y^{n_1} &  0          &\hdotsfor[2]{1}&   0    \\
    0      &0      & x^{m_2} &  x^{n_2}    &\hdotsfor[2]{1}&\vdots  \\
    0      &0      &\hdotsfor[2]{3}                        &   0    \\
    \vdots &\vdots &   0     &  0          &   x^{m_N}     & y^{n_N}\\
    0      &0      &   0     &  0          &    0          & x^{m_{N+1}}
   \end{pmatrix}.
  $$
  Hence, we can split $\widetilde{N}'(\boldsymbol{q}) = N(\boldsymbol{q})
  \oplus \Id_{R}$, which proves the claim.

  In the same way, we can show that the cokernel of
  $$N''(\boldsymbol{q}) := 
   \begin{pmatrix}
    y^{n_0}&   0     &  0          &\hdotsfor[2]{1}&   0        &0      \\
    x^{m_1}& y^{n_1} &  0          &\hdotsfor[2]{1}&   0        &0      \\
    0      & x^{m_2} &  y^{n_2}    &\hdotsfor[2]{1}&\vdots      &0      \\
    0      &\hdotsfor[2]{3}                        &   0        &0      \\
    \vdots &   0     &  0          &   x^{m_N}     & y^{n_N}    &\vdots \\
    0      &   0     &  0          &    0          & x^{m_{N+1}}&1      \\
    0      &\hdotsfor[2]{3}                        &   0        &1      \\
   \end{pmatrix}
  $$
  is isomorphic to $\mathcal{N}(\boldsymbol{q})$. In this sense,
  \begin{align*}
   \boldsymbol{q} &= n_{0}(m_{1},n_{1})\ldots(m_{N},n_{N})m_{N+1}\\
   \boldsymbol{q}' &= 0(0,n_{0})(m_{1},n_{1})\ldots(m_{N},n_{N})m_{N+1}\\
   \boldsymbol{q}'' &= n_{0}(m_{1},n_{1})\ldots(m_{N},n_{N})(m_{N+1},0)0
   \quad\text{ or }\\ 
   \boldsymbol{q}''' &= 0(0,n_{0})(m_{1},n_{1})\ldots(m_{N},n_{N})(m_{N+1},0)0.
  \end{align*}
  all define the same $R$-module $\mathcal{N}(\boldsymbol{q})$.
\end{remark}

Having these resolutions it is easy to compute the matrices of the 
multiplication with $X$ and $Y$. But it is much more instructive to visualise
these $R$-modules through a directed graph. We refer to them as \emph{band and
  string diagrams}. Such diagrams were introduced by Gelfand and Ponomarev
\cite{GelfandPonomarev}. They contain $\dim(V)$ vertices 
in the case of strings and $\frac{1}{m}\dim(V)$ vertices in the case of
bands. Each vertex corresponds to a subspace $V_{i}$ of $V$. In the case of
strings $\dim(V_{i})=1$, whereas in the case of bands one has $\dim(V_{i})=m$. 
The $V_{i}$ are quotients of subspaces of the direct summands of the middle
term in the resolution which defines the band or string module. 
These subspaces satisfy $X(V_{i})\subset V_{i-1}, Y(V_{i})\subset V_{i+1}$ 
and the restrictions of $X, Y$ on any $V_{i}$ are either isomorphisms or zero. 
If such a restriction is an isomorphism it is represented by an arrow in the
diagram. In the case of strings, one can choose a basis vector of each of the
$V_{i}$ such that all these isomorphisms map these basis vectors onto each
other. In this case, we label the arrows with $x$ respectively $y$ if it
represents the restriction of $X$ or $Y$, respectively. In the case of bands,
in addition, we have $Y(V_{n})\subset V_{1}$. After choosing appropriate bases
all but one of these isomorphisms can be represented by the identity matrix
$I_{m}$. Usually, we normalise the restriction of $Y$ to the space $V_{n}$
which corresponds to the ``last'' vertex to be the Jordan block
$J_{m}(\lambda)$. Equivalently, this could be done with any other non-zero
restriction of $Y$. Instead, we could normalise one of the non-zero
restrictions of $X$ to be $J_{m}(\lambda^{-1})$. In any case, the arrows in a
band diagram are labelled by $x$ or $y$ as above and the matrix which
represents the corresponding isomorphism. 

Details and proofs can be found in \cite{GelfandPonomarev}. 
We give here two examples.  

\begin{example}\label{ex:string}
  The string $\mathcal{N}(2(3,2)1)$ is represented by the diagram:

  \xymatrix@R-=12pt{
&&&&&{\boldsymbol{k}}  \ar@{->}[dl]_{x}\ar@{->}[dr]^{y}\\
  {\boldsymbol{k}}\ar@{->}[dr]^{y}&&&& {\boldsymbol{k}}\ar@{->}[dl]_{x}&& 
  {\boldsymbol{k}}\ar@{->}[dr]^{y}&&   {\boldsymbol{k}}\ar@{->}[dl]_{x} \\
& {\boldsymbol{k}}\ar@{->}[dr]^{y}&&   {\boldsymbol{k}}\ar@{->}[dl]_{x} &&&& 
  {\boldsymbol{k}} \\
&&{\boldsymbol{k}}
}
\end{example}

\begin{example}\label{ex:band}
  The band $\mathcal{M}(((2,2)(3,4)(1,3)),m,\lambda)$ is represented by:

\xymatrix@C-=4pt{
&&&&&&&{\boldsymbol{k}^{m}}\ar@{->}[dl]_{xI_{m}}\ar@{->}[dr]_{yI_{m}}\\
&&{\boldsymbol{k}^{m}}\ar@{->}[dl]_{xI_{m}}\ar@{->}[dr]_{yI_{m}}
&&&&{\boldsymbol{k}^{m}}\ar@{->}[dl]_{xI_{m}}
&&{\boldsymbol{k}^{m}}\ar@{->}[dr]_{yI_{m}}\\
&{\boldsymbol{k}^{m}}\ar@{->}[dl]_{xI_{m}}
&& {\boldsymbol{k}^{m}}\ar@{->}[dr]_{yI_{m}} 
&&{\boldsymbol{k}^{m}}\ar@{->}[dl]_{xI_{m}}
&&&& {\boldsymbol{k}^{m}}\ar@{->}[dr]_{yI_{m}}\\
 {\boldsymbol{k}^{m}}
&&&& {\boldsymbol{k}^{m}}
&&&&&&{\boldsymbol{k}^{m}}\ar@{->}[dr]_{yI_{m}}  
&&{\boldsymbol{k}^{m}}\ar@{->}[dl]_{xI_{m}}\ar@{->}[d]^{yI_{m}}\\
&&&&&&&&&&&{\boldsymbol{k}^{m}}
&{\boldsymbol{k}^{m}}\ar@{->}[d]^{yI_{m}}\\
&&&&&&&&&&&&{\boldsymbol{k}^{m}}\ar@{->}[uullllllllllll]_{yJ_{m}(\lambda)}
}
\end{example}

\section{Semi-Stable torsion free sheaves of degree zero}
\label{sec:semistable}

The aim of this section is twofold. First, we generalise results of \cite{FM}
on semi-stable vector bundles to semi-stable torsion free sheaves. The second
main achievement will be an explicit description of the functor $\mathbb{F}$
(see section \ref{sec:fm}) on semi-stable torsion free sheaves of degree
zero. We use the description of torsion sheaves supported at the singular point
which was given in section \ref{sec:tors}. 
Our main result is the following theorem.

\begin{theorem}\label{description}
  Let $\boldsymbol{E}$ be a nodal Weierstra{\ss} cubic and $E$ a semi-stable
  indecomposable torsion free sheaf of degree zero on $\boldsymbol{E}$ which
  is not of the form $\mathcal{L}\otimes\mathcal{F}_{m}$ with
  $\mathcal{L}\in\Pic(\boldsymbol{E})$.  
  \begin{itemize}
  \item[(a)] If $E$ is locally free, $E\cong
    \mathcal{B}(\boldsymbol{d},m,\lambda)$ with $n\ge 2$ and non-periodic
    $\boldsymbol{d}\in \mathbb{Z}^{n}$ of the form
    $$\boldsymbol{d}=(\overbrace{1,0,\ldots, 0}^{n_1},
                    \overbrace{-1, 0,\ldots,0}^{m_1},\ldots,
                    \overbrace{1,0,\ldots, 0}^{n_N},
                    \overbrace{-1, 0,\ldots,0}^{m_N}),$$
    where $N\ge 1, m_{i}\ge 1, n_{i}\ge 1$.
    Because $E$ is semi-stable, $\mathbb{F}(E)$ is a torsion sheaf. It is
    supported at $s\in\boldsymbol{E}$ and its stalk is the $R$-module
    $\mathcal{M}(\boldsymbol{q},m,(-1)^{n+N}\lambda)$, where
    $$\boldsymbol{q}=(n_{1},m_{1})(n_{2},m_{2})\ldots(n_{N},m_{N}).$$
  \item[(b)] If $E$ is not locally free, $E\cong \mathcal{S}(\boldsymbol{d})$
    with $n\ge 1$ and $\boldsymbol{d}\in\mathbb{Z}^{n}$ of  the form 
    $$\boldsymbol{d}=(\overbrace{0,\ldots 0}^{n_{1}}, -1, 
                      \overbrace{0,\ldots, 0}^{m_{1}}, 1,
                      \overbrace{0,\ldots, 0}^{n_{2}}, -1, \ldots, 1, 
                      \overbrace{0,\ldots, 0}^{n_{N}}, -1,
                      \overbrace{0,\ldots, 0}^{m_{N}}),$$
    where $N\ge 1, n_{i}\ge 0, m_{i}\ge 0$.
    Because $E$ is semi-stable, $\mathbb{F}(E)$ is a torsion sheaf. It is
    supported at $s\in\boldsymbol{E}$ and its stalk is the $R$-module
    $\mathcal{N}(\boldsymbol{q})$ with $$\boldsymbol{q} =
    0(n_{1},m_{1}+1)(n_{2}+1,m_{2}+1)\ldots(n_{N}+1,m_{N})0.$$
    If $n_{1}=0$ or $m_{N}=0$, we apply remark \ref{dummy} and cancel some
    zeroes at the end. If $N=1,2$ this might lead to $\boldsymbol{q} = n()m$.
  \end{itemize}
\end{theorem}

The indecomposable semi-stable torsion free sheaves of degree zero which are
excluded in the theorem were considered in corollary \ref{rankone}. The
sheaves $E$ which were considered in corollary \ref{rankone} all satisfy
$\ell(\mathbb{F}(E)_{s})\le 1$. The only sheaf which is treated in corollary 
\ref{rankone} and in theorem \ref{description} is the sheaf $\mathcal{S}(-1)$
which satisfies $\ell(\mathbb{F}(\mathcal{S}(-1))_{s})=1$. 

For strongly indecomposable vector bundles, this is the case $N=1$ of
theorem \ref{description} (a), the result was known to Friedman and Morgan, see
\cite{FM}, Cor. 2.3.2.

\begin{proof}[Proof of Theorem \ref{description}]
  Throughout the proof we use the notation introduced in section \ref{sec:tf}. 
  The idea of the proof is the following.
  For the sheaves $E$ which are described in the theorem, we calculate
  $\text{ev}\sphat_{s}$. This is the completion of the germ at the singular
  point $s\in\boldsymbol{E}$ of the evaluation map $\text{ev}:
  H^{0}(E(p_{0}))\otimes \mathcal{O} \rightarrow E(p_{0})$. It turns out that
  $\text{ev}\sphat_{s}$ is injective. This suffices to prove the injectivity
  of the evaluation map. Using $h^{0}(E(p_{0})) = \rk(E(p_{0}))$, this implies
  that $\coker(\text{ev})$ is a torsion sheaf. By remark \ref{computeF}, this
  cokernel is isomorphic to $\mathbb{F}(E)$. Finally, using corollary
  \ref{rankone}, theorem \ref{ss} and the results collected in section
  \ref{sec:tors}, we can prove the theorem. 

  From lemma \ref{cohomology},
  \begin{align*}
    \mathcal{B}(\boldsymbol{d},m,\lambda) \otimes \mathcal{O}(p_{0}) &\cong
    \mathcal{B}(\boldsymbol{d}+\boldsymbol{1},m,\lambda)\quad \text{ and }\\
    \mathcal{S}(\boldsymbol{d}) \otimes \mathcal{O}(p_{0})  &\cong
    \mathcal{S}(\boldsymbol{d}+\boldsymbol{1}),
  \end{align*}
  we obtain $H^{1}(E(p_{0}))=0$ and $h^{0}(E(p_{0}))=\rk(E(p_{0})) = mn$ for
  the sheaves $E$ described in the theorem. Here we let $m=1$ if
  $E=\mathcal{S}(\boldsymbol{d})$.

  Our first goal is the computation of $\ker(\text{ev}_{s})$ and
  $\coker(\text{ev}_{s})$. To identify these, it is sufficient to know its
  completion, because the completion of a Noetherian local ring is faithfully
  flat (\cite{Matsumura}, Thm. 8.14, see also lemma \ref{completion}).

  If $L$ is a line bundle on $\boldsymbol{E_{n}}$ and $\sigma\in H^{0}(L)$ a
  global section, the completion $\text{ev}\sphat_{s}$ at $s$ of the
  evaluation map $$\text{ev}: H^{0}(\pi_{n \ast}L)\otimes \mathcal{O}
  \rightarrow \pi_{n \ast}L$$ is described as follows.

  Because $\pi_{n}$ is \'etale, we have an isomorphism $(\pi_{n
  \ast}L)\sphat_{s} \cong \bigoplus_{\nu=1}^{n} L\sphat_{s_{\nu}}$. Using
  trivialisations of $L$ about each $s_{\nu}$, we obtain an isomorphism
  $(\pi_{n \ast}L)\sphat_{s} \cong \bigoplus_{\nu=1}^{n} R_{\nu}$. The
  $R$-module structure on it is given by (\ref{isom}) (see section
  \ref{etalecovers}), so that we identify it with $R^{n}$. The completed
  evaluation map sends the section $\sigma$ to the vector in $R^{n}$ whose
  $\nu$-th component is obtained by localising $\sigma$ at $s_{\nu}$.

  More explicitly, suppose $\sigma$ is represented by homogeneous polynomials
  $f_{\nu}(x_{\nu},y_{\nu})\in H^{0}(D_{\nu}, \mathcal{O}(d_{\nu}+1))$ which
  satisfy the gluing condition (\ref{gluing}) in section \ref{linebundles}. If
  $\nu\ne n$, the local description of $\sigma$ at the two preimages of the
  singularity $s_{\nu}$ is $(f_{\nu}(x_{\nu},1),f_{\nu+1}(1,y_{\nu+1})) \in
  \boldsymbol{k}[[x_{\nu}]] \oplus \boldsymbol{k}[[y_{\nu+1}]]$. The gluing
  condition ensures that these elements are indeed in $R_{\nu}$.
  The corresponding element in $R$, which is the $\nu$-th component of
  $\text{ev}\sphat_{s}(\sigma)$, is 
  $g_{\nu}=f_{\nu}(x,1)+f_{\nu+1}(1,y)-f_{\nu}(0,1)\in R$. At $s_{n}$, the
  gluing condition is $f_{n}(0,1)= \lambda f_{1}(1,0)$. The corresponding
  section on the normalisation is locally given by $(f_{n}(x_{n},1),\lambda
  f_{1}(1,y_{1})) \in \boldsymbol{k}[[x_{n}]] \oplus
  \boldsymbol{k}[[y_{1}]]$. The corresponding element in $R_{n}$ is represented
  by $g_{n}=f_{n}(x,1) +\lambda f_{1}(1,y)- f_{n}(0,1)\in R$.

  If $L$ is a line bundle on $\boldsymbol{I_{n}}$, the
  completion of
  $$\text{ev}: H^{0}(p_{n \ast}L)\otimes \mathcal{O} \rightarrow p_{n \ast}L$$
  at $s\in\boldsymbol{E}$ has a similar description. The main difference
  appears at $s_{0}$ and $s_{n}$. Again, $(p_{n \ast}L)\sphat_{s} \cong 
  \bigoplus_{\nu=0}^{n} R_{\nu}$, but now we have $R_{0}\cong
  \boldsymbol{k}[[y]]$, $R_{n} \cong \boldsymbol{k}[[x]]$ and all other
  $R_{\nu} \cong R$. The components $g_{1},\ldots,g_{n-1}$ are computed as
  above. But now, $g_{0} = f_{1}(1,y) \in R_{0}\cong \boldsymbol{k}[[y]]$ and
  $g_{n} = f_{n}(x,1) \in R_{n}\cong \boldsymbol{k}[[x]]$.

  Below, we use these descriptions to calculate a matrix representation of
  $\text{ev}\sphat_{s}$, which will be used to determine $\mathbb{F}(E)$.

  \noindent
  \underline{\textsl{The case $\mathcal{B}(\boldsymbol{d},1,\lambda)$.}} 

  Because $d_{\nu}\in\{-1,0,1\}$ for all $\nu$, it is easy to describe a basis
  of $H^{0}(E(p_{0}))$. Such a basis should have $n$ elements, because
  $m=1$. Observe that the vectors $\boldsymbol{d}\in\mathbb{Z}^{n}$ we study,
  satisfy $n\ge2, d_{1}=1$ and $d_{2}, d_{n}\in\{-1,0\}$. Below, we give all
  non-zero components $f_{\nu}$ for the basis elements we have chosen. If
  $n\ge 3$ we have:

  \smallskip
  \emph{type A.} Any $\nu$ with $d_{\nu}=1$ gives a basis vector with
  components: \label{basis}
  $$f_{\nu}=x_{\nu}y_{\nu}$$

  \smallskip
  \emph{type B.} Any $\nu$ with $d_{\nu}=-1$ gives a basis vector with
  components: 

  \begin{minipage}[c]{.4\textwidth}
  \begin{align*}
    \quad\underline{\nu\ne1,n}\quad
    f_{\nu-1} &= y_{\nu-1}^{d_{\nu-1}+1}\\
    f_{\nu}   &= 1\\
    f_{\nu+1} &= x_{\nu+1}^{d_{\nu+1}+1}
  \end{align*}
  \end{minipage}\hfill
  \begin{minipage}[c]{.4\textwidth}
  \begin{align*}
    \underline{\nu=n}\quad
    f_{n-1} &= y_{n-1}^{d_{n-1}+1}\\
    f_{n}   &= 1\\
    f_{1} &= \lambda^{-1} x_{1}^{2}
  \end{align*}
  \end{minipage}

  \smallskip
  \emph{type C.}
  Any $\nu$ with $d_{\nu}\ne -1 \ne d_{\nu+1}$ gives a basis vector: 

  \begin{minipage}[t]{.4\textwidth}
  \begin{align*}
    \underline{\nu\ne n}\quad\quad
    f_{\nu} &= y_{\nu}^{d_{\nu}+1}\\
    f_{\nu+1} &= x_{\nu+1}^{d_{\nu+1}+1}
  \end{align*}
  \end{minipage}\hfill
  \begin{minipage}[t]{.4\textwidth}
  \begin{align*}
    \underline{\nu=n}\quad\quad
    f_{n} &= y_{n}\\
    f_{1} &= \lambda^{-1} x_{1}^{2}
  \end{align*}
  \end{minipage}

  \bigskip
  It is easy to see that we obtained precisely $n=\rk(E)$ basis vectors. In
  the case $n=2$, there is no basis vector of \emph{type C} and the only
  change to be made occurs in \emph{type B}, where the basis vector is
  $(\lambda^{-1} x_{1}^{2}+y_{1}^{2}, 1)$.  

  As explained above, the completed evaluation map sends these sections to
  vectors in $R^{n}$ with components
  \begin{align*}
    g_{\nu} &= f_{\nu}(x,1)+f_{\nu+1}(1,y)-f_{\nu}(0,1) &\text{if }\; \nu<n,\\
    g_{n}   &= f_{n}(x,1) +\lambda f_{1}(1,y) - f_{n}(0,1).
  \end{align*}
  If $n\ge 4$ the basis above yields the following elements $(g_{1}, g_{2},
  \ldots, g_{n})\in R^{n}$. Again, we write down the non-zero components only: 

  \emph{type A:} ($d_{\nu}=1$)
  \begin{align*}
    (\nu&\ne1)&(\nu&=1)\\
    g_{\nu-1} &= y&g_{n}   &= \lambda y\\
    g_{\nu}   &=x&g_{1}   &=x 
  \end{align*}

  \emph{type B:} ($d_{\nu}=-1$)
  \begin{align*}
    (\nu&\ne1,2,n)&(\nu&=n)&(\nu&=2)\\
    g_{\nu-2} &= y^{d_{\nu-1}+1} &g_{n-2} &= y^{d_{n-1}+1}  
                                                &g_{n} &= \lambda y^{2}\\
    g_{\nu-1} &= 1               &g_{n-1} &= 1  &g_{1} &= 1\\
    g_{\nu}   &= 1               &g_{n}   &= 1  &g_{2}   &= 1\\
    g_{\nu+1} &= x^{d_{\nu+1}+1} &g_{1}   &= \lambda^{-1} x^{2}
                                                &g_{3}   &= x^{d_{3}+1}
  \end{align*}

  \emph{type C:} ($d_{\nu}\ne -1 \ne d_{\nu+1}$)
  \begin{align*}
    (\nu&\ne n,1)& (\nu&=n)&  (\nu&=1)\\
    g_{\nu-1} &= y^{d_{\nu}+1} &    g_{n-1} &= y & g_{n} &= \lambda y^{2}\\
    g_{\nu} &= 1               &    g_{n} &= 1   & g_{1} &= 1\\
   g_{\nu+1} &= x^{d_{\nu+1}+1}&    g_{1} &= \lambda^{-1} x^{2}&
                                            g_{2} &=  x
  \end{align*}

  The completed evaluation map as a mapping $R^{n}\rightarrow{ R^{n}}$ is
  described by the matrix whose columns are the vectors
  $(g_{1},\ldots,g_{n})$. Our computation below aims at showing that this
  mapping is injective and its cokernel is isomorphic to the cokernel of
  the matrix $M(\boldsymbol{q},1,(-1)^{n+N}\lambda)$. We might reduce our
  matrix using elementary operations of rows and columns (with coefficients
  from $R$). In addition, we might erase a row and a column if the entry which
  is in both of them is a unit in $R$ and all other entries in this row and
  this column are equal to zero. This operation corresponds to splitting an
  isomorphism $R \rightarrow R$ as a direct summand from $R^{n}\rightarrow{
  R^{n}}$. We call two matrices \emph{equivalent} if they have isomorphic
  kernels and cokernels.

  Observe that any vector $g$ of \emph{type A} has the
  property that $x\cdot g$ has precisely one non-zero component and this is
  equal to $x^{2}$. Similarly, $y\cdot g$ has $y^{2}$ (or $\lambda y^{2}$, if
  $\nu=1$) as its only non-zero component. The position of the value $x^{2}$
  in $x\cdot g$ is the component 
  with number $\nu$ precisely if $d_{\nu}=1$. The value $y^{2}$ sits in $y\cdot
  g$ at the component with number $\nu$ precisely if $d_{\nu+1}=1$.
  On the other hand, values $x^{2}$, $\lambda x^{2}$, $y^{2}$ or $\lambda^{-1}
  y^{2}$ occur as components $g_{\mu}$ of vectors of \emph{type B} and
  \emph{type C} precisely at the same positions. Hence, we can annihilate in
  our matrix all entries which involve $x^{2}$ or $y^{2}$.

  If $n=2$ we have $\boldsymbol{d}=(1,-1)$ and the two vectors are $(x,\lambda
  y)$ and $(\lambda^{-1}x^{2}+1, \lambda y^{2}+1)$. If $n=3$ and
  $\boldsymbol{d}=(1,0,-1)$, the three image vectors are $(x,0,\lambda y),
  (\lambda^{-1}x^{2}+y,1,1), (1, x, \lambda y^{2})$. If $\boldsymbol{d} = 
  (1,-1,0)$, the three image vectors are $(x,0,\lambda y), (1,1,x+\lambda
  y^{2}), (\lambda^{-1}x^{2}, y, 1)$. These three cases are dealt with easily
  by hand and are left to the reader. Suppose $n\ge 4$ for the rest of the
  proof. 

  Each part of $\boldsymbol{d}$ of the form
  $\overbrace{1,0,\ldots,0}^{n_{k}},\overbrace{-1,0,\ldots,0}^{m_{k}},1$ gives
  rise to the following part of our matrix:

  \begin{center}\label{matrixdia}\begin{tabular}[t]{c|ccccccccccccc|c}
   $1$&\textcircled{$y$}&$\underline{0}$&&&&&&&&&&&&\\ \hline\hline
   $\underline{0}$&$x$&$1$&$y$&&&&&&&&&&&\\
   &&$x$&$1$&$y$&&&&&&&&&&\\
   &&&$\cdot$&$\cdot$&$\cdot$&&&&&&&&&\\
   &&&&$\cdot$&$\cdot$&$\cdot$&&&&&&&&\\
   &&&&&$x$&$1$&$y$&&&&&&&\\
   &&&&&&$x$&$1$&$0$&&&&&&\\ \hline
   &&&&&&$0$&$1$&$y$&&&&&&\\
   &&&&&&&$x$&$1$&$y$&&&&&\\
   &&&&&&&&$\cdot$&$\cdot$&$\cdot$&&&&\\
   &&&&&&&&&$\cdot$&$\cdot$&$\cdot$&&&\\
   &&&&&&&&&&$x$&$1$&$y$&&\\
   &&&&&&&&&&&$x$&$1$&\textcircled{$y$}&$\underline{0}$\\ \hline\hline
   &&&&&&&&&&&&$\underline{0}$&$x$&$1$\\
  \end{tabular}\end{center}
  If $N=1$, this has to be changed slightly. In this case the matrix consists
  just of the part inside the large box. Moreover, the last column inside the
  box has to be erased and its content must be added to the first column. In
  this way, the encircled entry $y$ in the last row appears in the first
  column. With these changes, the procedure described below applies to the
  case $N=1$ as well. 
  
  Because the rows of our matrix correspond to singular points, the space
  between two rows corresponds to a component of $\boldsymbol{E_{n}}$ and so to
  one of the components of $\boldsymbol{d}$. In this sense, the double lines
  correspond to $d_{\nu}=1$ and the single horizontal line in the middle
  corresponds to $d_{\nu}=-1$. The number of rows inside the upper part of the
  box is $n_{k}$. If $n_{k}=1$ only the last row of the upper part is
  present. Similarly, if $m_{k}=1$ the only row which is present in the lower
  part of the depicted portion of the matrix is the first row of the lower
  part. The lower part contains $m_{k}$ rows. The underlined zero entries
  indicate the positions where we deleted (multiples of) $x^{2}$ or $y^{2}$. If
  $k=N$, the encircled entry $y$ in the last row above the lower double line
  is to be replaced by $\lambda y$. If $k=1$, the encircled entry $y$ in
  first row has to be replaced by $\lambda y$, but this does not influence our
  reduction of the part inside the boxes. All non-zero entries are visible, if
  they belong to a row or column, a part of which is contained inside the
  depicted box. Outside this box, these are just the entry $y$ in the first
  row and the entry $x$ in the last row.

  We use that we have $xy=0$ in $R$ and subtract successively $x$ or $y$ times
  a column from one of its neighbours. In this way, we find that our original
  matrix is equivalent to one where the portion above is replaced by

  \begin{center}\begin{tabular}[t]{c|ccccccccccccc|c}
   $1$&\textcircled{$y$}&$\underline{0}$&&&&&&&&&&&&\\ \hline\hline
   $\underline{0}$&$0$&$1$&$0$&&&&&&&&&&&\\
   &&$0$&$1$&$0$&&&&&&&&&&\\
   &&&$\cdot$&$\cdot$&$\cdot$&&&&&&&&&\\
   &&&&$\cdot$&$\cdot$&$\cdot$&&&&&&&&\\
   &&&&&$0$&$1$&$0$&&&&&&&\\
   &$-(-x)^{n_{k}}$&&&&&$0$&$1$&$0$&&&&&&\\ \hline
   &&&&&&$0$&$1$&$0$&&&&&$-\lambda(-y)^{m_{k}}$&\\
   &&&&&&&$0$&$1$&$0$&&&&&\\
   &&&&&&&&$\cdot$&$\cdot$&$\cdot$&&&&\\
   &&&&&&&&&$\cdot$&$\cdot$&$\cdot$&&&\\
   &&&&&&&&&&$0$&$1$&$0$&&\\
   &&&&&&&&&&&$0$&$1$&$0$&$\underline{0}$\\ \hline\hline
   &&&&&&&&&&&&$\underline{0}$&$x$&$1$\\
  \end{tabular}\end{center}

  Of course, here and below, $\lambda$ must be replaced by $1$ if $k\ne N$. If
  $k=1$, the encircled entry $y$ in the first row should be replaced by
  $\lambda y$. We can now erase the $n_{k}-1$ rows and columns which meet at
  the unit matrix in the upper part and similarly $m_{k}-1$ rows from the lower
  part with their corresponding columns. We obtain an equivalent matrix, if we
  replace the original portion of our matrix by:
  \begin{center}\begin{tabular}[t]{c|ccc|c}
   $1$&$y$&&&\\ \hline\hline
   &$-(-x)^{n_{k}}$&$1$&$0$&\\ \hline
   &$0$&$1$&$-\lambda(-y)^{m_{k}}$&\\ \hline\hline
   &&&$x$&$1$\\
  \end{tabular}\end{center}
  In a final step we reduce this to
  \begin{center}\begin{tabular}[t]{c|cc|c}
   $1$&$y$&&\\ \hline\hline
   &$(-x)^{n_{k}}$&$-\lambda(-y)^{m_{k}}$&\\ \hline\hline
   &&$x$&$1$\\
  \end{tabular}\end{center}
  If $N=1$, there is just one column left and the entry is
  $(-x)^{n_{1}}-\lambda(-y)^{m_{1}}$ or, equivalently,
  $x^{n_{1}}+(-1)^{n_{1}+m_{1}+1}\lambda y^{m_{1}}$.

  If $N>1$, this process does not change any non-zero entry outside the area
  between the two double lines. Moreover, we can perform the same procedure
  even if the four non-zero entries $1,y,x,1$ in the picture, which are outside
  the box in which the changes take place, are replaced by other
  values. Therefore, each of the $N$ blocks as described above can be replaced
  by the corresponding $1\times 2$ matrix which was obtained at the end of the
  procedure. Up to this point we think of the numbers of the rows of this
  matrix as cyclic subscripts. In other words, we do not fix a preference which
  subscript is chosen as the first row in our matrix. But now, at the end,
  when we write down the complete matrix, we chose to write the row which
  corresponds to $s_{1}$ as our first row. In this way, we obtain that the
  original matrix is equivalent to the following matrix of size $N\times N$:
  $$\begin{pmatrix}
    (-x)^{n_{1}}             &-(-y)^{m_{1}}&0            &&&\\
    0                        &(-x)^{n_{2}} &-(-y)^{m_{2}}&&&\\
    0                        &0            &(-x)^{n_{3}} &&&\\
                             &             &             &\ddots&&\\
                             &             &&&(-x)^{n_{N-1}} &-(-y)^{m_{N-1}}\\
    -\lambda(-y)^{m_{N}}&             &&&0              &(-x)^{n_{N}} 
  \end{pmatrix}.$$
  Using $n=\sum_{\nu=1}^{N} (m_{\nu}+n_{\nu})$ this is easily seen to be
  equivalent to
  $$\begin{pmatrix}
    x^{n_{1}}                  &y^{m_{1}} &0             &&&\\
    0                          &x^{n_{2}} &y^{m_{2}}     &&&\\
    0                          &0         &x^{n_{3}}     &&&\\
                               &          &              &\ddots&&\\
                               &          & &&x^{n_{N-1}}&y^{m_{N-1}}\\
    (-1)^{n+N}\lambda y^{m_{N}}&          & &&0          &x^{n_{N}} 
  \end{pmatrix},$$
  which is $M(\boldsymbol{q},1,(-1)^{n+N}\lambda)$. Hence, we have shown
  $$\coker(\text{ev}_{s})\sphat \cong
  \mathcal{M}(\boldsymbol{q},1,(-1)^{n+N}\lambda) \quad\text{ and }\quad
  \ker(\text{ev}_{s})\sphat = 0.$$
  This implies $\ker(\text{ev}_{s}) = 0$ and, because $\boldsymbol{E}$ is
  irreducible, $\ker(\text{ev})$ has rank zero. But, as a subsheaf of a
  torsion free sheaf, $\ker(\text{ev})$ is torsion free itself, hence
  zero. Remark \ref{computeF} implies now
  $\mathbb{F}(\mathcal{B}(\boldsymbol{d},1,\lambda)) \cong \coker(\text{ev})$.
  Because $h^{0}(E(p_{0})) = \rk(E(p_{0}))$, this cokernel is a torsion sheaf
  and, with $\boldsymbol{d}$ and $\boldsymbol{q}$ as described in the theorem,
  we have 
  $$\mathbb{F}(\mathcal{B}(\boldsymbol{d},1,\lambda))\sphat_{s} \cong
  \mathcal{M}(\boldsymbol{q},1,(-1)^{n+N}\lambda),$$ which is an
  indecomposable module of finite length. Theorem \ref{ss} implies now that
  $\mathcal{B}(\boldsymbol{d},1,\lambda)$ is semi-stable, if $\boldsymbol{d}$
  is of the form which is described in the theorem. To see that
  $\mathcal{B}(\boldsymbol{d},1,\lambda)$ is indecomposable, we have to
  exclude that the support of
  $\mathbb{F}(\mathcal{B}(\boldsymbol{d},1,\lambda))$ contains a regular point
  $p\in\boldsymbol{E}$. To see this, it is sufficient to show
  $$\Hom(\boldsymbol{k}(p),\mathbb{F}(\mathcal{B}(\boldsymbol{d},1,\lambda))) =
  0.$$
  By remark \ref{braid}, $\mathbb{F}$ is an equivalence of categories and
  by corollary \ref{rankone}, $\boldsymbol{k}(p) \cong \mathbb{F}(L)$, where
  $L$ is a line bundle of degree zero on $\boldsymbol{E}$. Thus,
  \begin{align*}    
  \Hom(\boldsymbol{k}(p),\mathbb{F}(\mathcal{B}(\boldsymbol{d},1,\lambda)))
  &\cong \Hom(L, \mathcal{B}(\boldsymbol{d},1,\lambda)) \cong\\
  \Hom(L, \pi_{n\ast}(\mathcal{L}(\boldsymbol{d},\lambda))) &\cong
  \Hom(\pi_{n}^{\ast}(L), \mathcal{L}(\boldsymbol{d},\lambda)).
  \end{align*}
  As $L$ has degree zero, $\pi_{n}^{\ast}(L)$ is of degree zero on each
  component of $\boldsymbol{E_{n}}$. Hence,
  $\Hom(\pi_{n}^{\ast}(L), \mathcal{L}(\boldsymbol{d},\lambda)) \cong
  H^{0}(\mathcal{L}(\boldsymbol{d},\lambda'))$, with some
  $\lambda'\in\boldsymbol{k}^{\times}$. But
  $H^{0}(\mathcal{L}(\boldsymbol{d},\lambda')) = 0$ (same proof as lemma
  \ref{cohomology}), hence $\mathbb{F}(\mathcal{B}(\boldsymbol{d},1,\lambda))$
  is supported at $s$ only. This shows that we identified all indecomposable
  semi-stable torsion free sheaves of degree zero whose Fourier-Mukai image is
  one of the band modules $\mathcal{M}(\boldsymbol{q},1,\lambda)$.

  \noindent
  \underline{\textsl{The case $\mathcal{B}(\boldsymbol{d},m,\lambda)$.}} 

  Using induction on $m$ and the exact sequences
  $$0  \rightarrow \mathcal{B}(\boldsymbol{d},1,\lambda)
       \rightarrow \mathcal{B}(\boldsymbol{d},m+1,\lambda)
       \rightarrow \mathcal{B}(\boldsymbol{d},m,\lambda)
       \rightarrow 0,$$
  we obtain the semi-stability of $\mathcal{B}(\boldsymbol{d},m,\lambda)$ for
  any $m\ge 1$, provided $\boldsymbol{d}$ is one of the vectors appearing in
  the theorem. Theorem \ref{ss} implies that
  $\mathbb{F}(\mathcal{B}(\boldsymbol{d},m,\lambda))$ is isomorphic to the
  cokernel of the evaluation map. Using the same exact sequences, exactness of
  $\mathbb{F}$ and induction on $m$ imply that the support of
  $\mathbb{F}(\mathcal{B}(\boldsymbol{d},m,\lambda))$ is the singular point
  $s$.
  
  If $m>1$ we have $\mathcal{B}(\boldsymbol{d},m,\lambda) \cong
  \pi_{n\ast}(\mathcal{L}(\boldsymbol{d},\lambda) \otimes
  \pi_{n}^{\ast}\mathcal{F}_{m})$. The vector bundle
  $\mathcal{L}(\boldsymbol{d},\lambda) \otimes \pi_{n}^{\ast}\mathcal{F}_{m}$
  on $\boldsymbol{E_{n}}$ is obtained by gluing the bundles
  $\mathcal{O}(d_{\nu})^{\oplus m}$ in a similar way as we obtained the
  bundles $\mathcal{L}(\boldsymbol{d},\lambda)$: with the same convention and
  coordinates as before, over $s_{1},\ldots,s_{n-1}$ we glue with the identity
  matrix $I_{m}$ and over $s_{n}$ we glue with the matrix
  $J_{m}(\lambda)$. Actually, the definition of the bundle implies that the
  gluing matrix over $s_{n}$ should be $\lambda J_{m}(1)^{n}$, but this matrix
  is similar to $J_{m}(\lambda)$. (The assumption $\cha(\boldsymbol{k}) = 0$
  is essential here!) Sections of
  $\mathcal{B}(\boldsymbol{d},m,\lambda)$ are now given by $n$-tuples 
  $(f_{1},\ldots,f_{n}) \in \oplus H^{0}(\mathcal{O}(d_{\nu}))^{\oplus m}$
  which satisfy the gluing condition (\ref{gluing}) with $\lambda$ replaced by
  $J_{m}(\lambda)$ in the second equation.

  As a result, we obtain a basis of
  $H^{0}(\mathcal{B}(\boldsymbol{d}+\boldsymbol{1},m,\lambda))$ if we
  interpret the description of the basis of
  $H^{0}(\mathcal{B}(\boldsymbol{d}+\boldsymbol{1},1,\lambda))$ 
  on page \pageref{basis} in the following way: each basis vector in the case
  $m=1$ gives rise to $m$ basis vectors by replacing $x^{k},y^{k},\lambda
  y^{k}$ by $x^{k}I_{m}, y^{k}I_{m}$ and $y^{k}J_{m}(\lambda)$
  respectively. The $m$ columns produced this way are the new basis vectors.

  We can now carry out the same procedure as in the case $m=1$. The only
  differences are now that the entries in our matrix are blocks of size
  $m\times m$ and that the factor $\lambda$ has to be replaced by the matrix
  $J_{m}(\lambda)$. This does not influence any of the reduction steps, so
  that we finally obtain the matrix
  $M(\boldsymbol{q},m,(-1)^{n+N}\lambda)$. This proves
  $$\mathbb{F}(\mathcal{B}(\boldsymbol{d},m,\lambda))\sphat_{s} \cong 
  \mathcal{M}(\boldsymbol{q},m,(-1)^{n+N}\lambda).$$ 
  As seen above, for any $m\ge 1$, the support of the 
  sheaf $\mathbb{F}(\mathcal{B}(\boldsymbol{d},m,\lambda))$ is the point $s$. 
  Using the results from section \ref{sec:tors} this shows
  that any indecomposable $R$-module $\mathcal{M}(\boldsymbol{q},m,\lambda)$ is
  isomorphic to  a module $\mathbb{F}(E)$ with
  $E\cong\mathcal{B}(\boldsymbol{d},m,(-1)^{n+N}\lambda)$ as given in the
  theorem.

  \noindent
  \underline{\textsl{The case $\mathcal{S}(\boldsymbol{d})$.}}

  If $E=\mathcal{S}(\boldsymbol{d})$, the computations are very similar to the
  computations in the case $E=\mathcal{B}(\boldsymbol{d},1,\lambda)$.
  We confine ourselves to describe the differences.
  First of all, an $n$-tuple
  $(f_{1}, f_{2}, \ldots, f_{n}) \in \bigoplus_{\nu=1}^{n} H^{0}(D_{\nu},
  \mathcal{O}(d_{\nu}+1))$ is an element of $H^{0}(E(p_{0}))$ if and only if 
  \begin{align*}
    f_{\nu}(0:1) &= f_{\nu+1}(1:0) &1\le \nu \le n-1.
  \end{align*}
  We start with an explicit description of a basis of $H^{0}(E(p_{0}))$. As
  before, only the non-zero entries of the vectors $(f_{1},\ldots,f_{n})$ are
  written down. Observe, $d_{1}\ne 1, d_{n}\ne 1, n\ge 2$. 

  \emph{type A.} Any $\nu$ with $d_{\nu}=1$ gives a basis vector with
  components: 
  \begin{align*}
    f_{\nu}   &=x_{\nu}y_{\nu}
  \end{align*}

  \emph{type B:} Any $\nu$ with $d_{\nu}=-1$ gives a basis vector with
  components: 
  \begin{align*}
    (\nu&\ne1,n)&(\nu&=n)&(\nu&=1)\\
    f_{\nu-1} &= y_{\nu-1}^{d_{\nu-1}+1} &f_{n-1} &=y_{\nu-1}^{d_{n-1}+1}&&\\
    f_{\nu}   &= 1                       &f_{n}   &= 1  &f_{1}&= 1\\
    f_{\nu+1} &= x_{\nu+1}^{d_{\nu+1}+1} &        &     &f_{2}&=x_{2}^{d_{2}+1}
  \end{align*}

  \emph{type C:} Any $\nu$ with $d_{\nu}\ne -1 \ne d_{\nu+1}$ gives a basis
  vector: 
  \begin{align*}
    (\nu&\ne n)& (\nu&=n)\\
    f_{\nu}   &= y_{\nu}^{d_{\nu}+1}&    f_{n} &= y_{n}   \\
    f_{\nu+1} &= x_{\nu+1}^{d_{\nu+1}+1}    &    
  \end{align*}

  \emph{type D:} If $d_{1}=0$ there is a basis vector with:
  \begin{align*}
    f_{1}   &= x_{1}
  \end{align*}

  Again, the number of these vectors is equal to $n$ and their images
  $g=(g_{0},\ldots,g_{n})$ under the completed evaluation map are computed by
  the procedure explained above. 

  We shall study parts of the matrix which represents the completed evaluation
  map $R^{n} \rightarrow \boldsymbol{k}[[y]] \oplus R^{n-1} \oplus
  \boldsymbol{k}[[x]]$.  As before, we look at those parts which correspond to
  portions of $\boldsymbol{d}$ the form
  $$1,\overbrace{0,\ldots,0}^{n_{k}}, -1, \overbrace{0,\ldots,0}^{m_{k}},1$$
  with $1\le k\le N$. If $k=1$ there is no leading $1$, whereas the trailing
  $1$ is missing in case $k=N$. If $1\ne k\ne N$, the corresponding portion
  of the matrix looks exactly like before. The number of rows in the upper
  part inside the box is now equal to $n_{k}+1$ and in the lower part equal to
  $m_{k}+1$. The same reduction steps as in the locally free case lead to
  
  \begin{center}\begin{tabular}[t]{c|cc|c}
   $1$&$y$&&\\ \hline\hline
   &$(-x)^{n_{k}+1}$&$-(-y)^{m_{k}+1}$&\\ \hline\hline
   &&$x$&$1$\\
  \end{tabular}\end{center}

  However, if $k=1$ the row above the upper double line is not present
  in the portion of the matrix depicted on page
  \pageref{matrixdia}. Furthermore, the first column inside the box is to be
  deleted as well. So, the upper part of the box contains $n_{1}+1$ rows and
  the upper left corner looks like: 
  \begin{center}\begin{tabular}[t]{|cccc}\hline
   $1$&$y$&$0$&$\cdots$\\ 
   $x$&$1$&$y$&\\ 
   $0$&$x$&$1$&$y$
  \end{tabular}\end{center}
  and (even if $n_{1}=0$) we reduce it to 
  \begin{center}\begin{tabular}[t]{|cc|c}\hline
   $1$&$0$&\\ 
   $-(-x)^{n_{1}}$&$-(-y)^{m_{1}+1}$&\\ \hline\hline
   $0$&$x$&$1$
  \end{tabular}\end{center}
  Similarly, if $k=N$, we have $n_{N}+1$ rows in the upper part and
  $m_{N}+1$ rows in the lower part and deal with a lower right corner or the
  form  
  \begin{center}\begin{tabular}[t]{cccc|}
   $x$&$1$&$y$&$0$\\
   &$x$&$1$&$y$\\ 
   $\cdots$&$0$&$x$&$1$\\ \hline
  \end{tabular}\end{center}
  and (even if $m_{N}=0$) the reduction leads to
  \begin{center}\begin{tabular}[t]{c|cc|}
   $1$&$y$&$0$\\ \hline\hline
   &$(-x)^{n_{N}+1}$&$-(-y)^{m_{N}}$\\ 
   &$0$&$1$\\ \hline
  \end{tabular}\end{center}
  We need to keep in mind here that the entry in the first row is an element of
  $\boldsymbol{k}[[y]]$ with trivial multiplication by $x\in R$ and the entry
  in the last row is in $\boldsymbol{k}[[x]]$ with trivial multiplication by
  $y\in R$. Multiplying rows and columns by appropriate powers of $(-1)$, we
  arrive at the following matrix, which is equivalent to the completed
  evaluation map: 
  $$\begin{pmatrix}
    1          &0          &\hdotsfor[2]{3}&0\\
    x^{n_{1}}  &y^{m_{1}+1}&0              &\hdotsfor[2]{2}&0\\
    0          &x^{n_{2}+1}&y^{m_{2}+1}    &0&\hdotsfor[2]{1}&0\\
    0          &0          &x^{n_{3}+1}    &\ddots&&\vdots\\
    \vdots     &           &               &\ddots&y^{m_{N-1}+1}&0\\
    0          &\hdotsfor[2]{2}            &0     &x^{n_{N}+1} &y^{m_{N}}\\
    0          &\hdotsfor[2]{3}                   &0           &1 
  \end{pmatrix}$$
  This matrix defines a mapping $R^{N+1}\rightarrow \boldsymbol{k}[[y]] \oplus
  R^{N} \oplus \boldsymbol{k}[[x]]$.

  If $n_{1}>0$ and $m_{N}>0$ this is precisely the matrix
  $N(\boldsymbol{q})$. Whereas, if $n_{1}=0$ or $m_{N}=0$, we obtained one of
  the matrices $N'(\boldsymbol{q})$ or $N''(\boldsymbol{q})$ of remark
  \ref{dummy}. This shows that the cokernel of this matrix is isomorphic to the
  module $\mathcal{N}(\boldsymbol{q})$, with $\boldsymbol{q}$ being the
  sequence described in the theorem. In addition, we see that
  $\text{ev}\sphat_{s}$ is injective.

  As before, we conclude $\mathbb{F}(\mathcal{S}(\boldsymbol{d})) \cong
  \coker(\text{ev})$. Our computation implies now 
  $\mathbb{F}(\mathcal{S}(\boldsymbol{d}))\sphat_{s} \cong
  \mathcal{N}(\boldsymbol{q})$ for $\boldsymbol{d}$ and $\boldsymbol{q}$ as in
  the theorem. Together with the results from section \ref{sec:tors} this shows
  that any indecomposable $R$-module $\mathcal{N}(\boldsymbol{q})$ has the form
  $\mathbb{F}(\mathcal{S}(\boldsymbol{d}))\sphat_{s}$ with 
  $\boldsymbol{d}$ and $\boldsymbol{q}$ as given in the theorem.
  To see that $\mathbb{F}(\mathcal{S}(\boldsymbol{d}))$ and, hence,
  $\mathcal{S}(\boldsymbol{d})$ are indecomposable, we proceed as before in
  the case $\mathcal{B}(\boldsymbol{d},1,\lambda)$. We have to verify
  $\Hom(p_{n}^{\ast}(L), \mathcal{L}(\boldsymbol{d})) = 0$ for any line bundle
  $L$ of degree zero on $\boldsymbol{E}$. But $p_{n}^{\ast}(L)
  \cong\mathcal{O}$  on $\boldsymbol{I_{n}}$. Hence, $\Hom(p_{n}^{\ast}(L),
  \mathcal{L}(\boldsymbol{d})) \cong H^{0}(\mathcal{L}(\boldsymbol{d}))
  =0$. This proves the indecomposability of $\mathcal{S}(\boldsymbol{d})$ and
  that $\mathbb{F}(\mathcal{S}(\boldsymbol{d}))$ is supported at $s$ only.

  The conclusion now is that any coherent indecomposable torsion module on
  $\boldsymbol{E}$ with support in $s$ has the form $\mathbb{F}(E)$ where $E$
  is a sheaf as described in the theorem. On the other hand, any
  indecomposable torsion sheaf of finite length which is supported at a
  regular point $x\in\boldsymbol{E}$ is isomorphic to
  $\mathcal{O}_{\boldsymbol{E},x}/\mathfrak{m}_{x}^{m}$ with $m\ge 1$ being
  the length of the sheaf. We have shown in corollary \ref{rankone} that these
  sheaves are isomorphic to $\mathbb{F}(\mathcal{L}(0,\lambda) \otimes
  \mathcal{F}_{m})$. Using theorem \ref{ss} this implies that we described all
  indecomposable semi-stable torsion free sheaves of rank zero.
\end{proof}

\begin{remark}
  Indecomposable semi-stable vector bundles of degree zero were
  characterised in \cite{FM}. We do not use their result here and produce a
  new proof of the fact that the vectors $\boldsymbol{d}$ as   described in
  the theorem correspond precisely to the semi-stable vector bundles. If $E$
  is not locally free, semi-stability of these sheaves was not known before. 
\end{remark}

\begin{remark}
  To avoid any possible confusion or misinterpretation, the structure of the
  vectors $\boldsymbol{d}\in\mathbb{Z}^{n}$ in both parts of theorem
  \ref{description} is the following: if all the zero components are removed,
  we obtain an alternating sequence of numbers $-1$ and $1$.
\end{remark}

\begin{remark}
  Although the description of the sequences $\boldsymbol{q}$ in part (b) of
  theorem \ref{description} seems to be a bit awkward, the description of the
  string diagram of $\mathbb{F}(\mathcal{S}(\boldsymbol{d}))$ in terms of
  $\boldsymbol{d}$ is straightforward: if $\boldsymbol{d}$ is of the form
  described in part (b) of theorem \ref{description}, the string diagram
  contains a vertex for each component $d_{\nu}$ of $\boldsymbol{d}$, the
  arrows connect neighbours only and form sequences pointing in the same
  direction which start at vertices corresponding to $d_{\nu}=-1$ and end at
  vertices which correspond to $d_{\nu}=1$. This applies to the sheaf
  $\mathcal{S}(-1)$ as well. In this case, the diagram consists of one vertex
  and no arrows.

  Literally the same description applies to give a description of the band
  diagram for $\mathbb{F}(\mathcal{B}(\boldsymbol{d},m,\lambda))$. Here we
  require $\boldsymbol{d} \ne \boldsymbol{0}$ to obtain a sheaf supported at
  the singularity. The only difference to the string case is that the
  components of $\boldsymbol{d}$ are considered to be in cyclic order, so that
  $d_{n}$ and $d_{1}$ are neighbours, hence connected by an arrow.
\end{remark}

\section{Matlis Duality}
\label{sec:matlis}

In this section we study the relationship between $\mathbb{F}(E)$ and
$\mathbb{F}(E^{\vee})$, where $E$ is a semi-stable torsion free sheaf of degree
zero on a Weierstra{\ss} cubic $\boldsymbol{E}$. The answer (theorem
\ref{dual}) involves Matlis duality.

We briefly recall the features of Matlis duality which will be used
later. Details and a more comprehensive treatment can be found in
\cite{Matsumura} and \cite{BrunsHerzog}. 

First, recall that a local Noetherian ring $(A,\mathfrak{m},\boldsymbol{k})$ is
called \emph{Gorenstein}, if its injective dimension is finite. Further, the
injective hull of an $A$-module $M$ is an injective $A$-module $E(M)$ which
contains $M$ such that $M\subset E(M)$ is an essential extension. To be an
essential extension means here that for any non-zero submodule $N\subset E(M)$
one has $M\cap N\ne 0$. 
The indecomposable injective $A$-modules are precisely the modules
$E(A/\mathfrak{p})$, where $\mathfrak{p}\subset A$ is a prime ideal
\cite{Matsumura}, Thm. 18.4.  

\begin{lemma}\label{injective}
  Let $(A,\mathfrak{m},\boldsymbol{k})$ be a local Gorenstein ring of dimension
  one. By $K:=\Quot(A)$ we denote its total quotient ring. Let
  $\mathfrak{p}_{1}, \ldots, \mathfrak{p}_{m}$ be the minimal prime ideals of
  $A$. 
  \begin{itemize}
  \item[(i)] One has $K/A\cong E(\boldsymbol{k})$ and $K \cong E(A) \cong 
    \bigoplus_{\nu=1}^{m} E(A/\mathfrak{p}_{\nu})$.
  \item[(ii)] If $M$ is an $A$-module of finite length (or, equivalently, a
    torsion module, i.e. $M\otimes_{A} K =0$), then
    $$\Hom_{A}(M,K/A) \cong \Ext_{A}^{1}(M,A).$$
  \end{itemize}
\end{lemma}

\begin{proof}
  (i) This follows easily from the standard theory, see
  \cite{Matsumura}, Ch. 18. 

  (ii) The assumption that $M$ is torsion implies
  $\Hom_{A}(M,E(A/\mathfrak{p}_{\nu}) )=0$ for all $\nu$. 
  Hence, applying the functor $\Hom_{A}(M,-)$ to the minimal injective
  resolution of $A$ (see \cite{Matsumura}, Thm. 18.8)
  $$0 \rightarrow A \rightarrow \bigoplus_{\nu=1}^{m} E(A/\mathfrak{p}_{\nu}) 
                   \rightarrow E(\boldsymbol{k}) \rightarrow 0,$$
  yields the desired isomorphism.
\end{proof}

\begin{definition}
  If $(A,\mathfrak{m},\boldsymbol{k})$ is a local Noetherian ring, the functor
  $$\mathbb{M}_{A}(-) := \Hom_{A}(-,E(\boldsymbol{k}))$$ is called the
  \emph{Matlis functor}. 
\end{definition}

The main result of Matlis (see \cite{Matsumura}, Thm. 18.6) states:
\begin{itemize}
\item[(a)] For any $A$-module $M$, the canonical map
  $$M\rightarrow \mathbb{M}_{A}(\mathbb{M}_{A}(M))$$ is injective. 
\item[(b)] If $M$ is an $A$-module of finite length $\ell(M)$, one has
  $\ell(\mathbb{M}_{A}(M)) = \ell(M)$ and the canonical map in (a) is an
  isomorphism. 
\item[(c)] If $A$ is complete, $\mathbb{M}_{A}$ is an anti-equivalence between
  the categories of Artinian and Noetherian $A$-modules. For both types of
  modules, the canonical map in (a) is an isomorphism. 
\end{itemize}
These results can be used to prove Grothendieck's local duality theorem,
see \cite{BrunsHerzog}, Sect. 3.5. 

\begin{lemma}\label{kdual}
  Let $(A,\mathfrak{m},\boldsymbol{k})$ be a local Gorenstein ring of
  dimension one. Suppose $E$ is an injective $A$-module which satisfies
  $\dim_{\boldsymbol{k}} \Hom_{A}(\boldsymbol{k},E) =1$. Then there is an
  isomorphism of functors on the category of $A$-modules of finite length:
  $$\mathbb{M}_{A}(-) \cong \Hom_{A}(-,E).$$ 
\end{lemma}

\begin{proof}
  Using \cite{Matsumura}, Thm.\/ 18.5, our assumption implies $E\cong
  E(\boldsymbol{k}) \oplus E'$, where $E'$ is a direct sum of injective modules
  $E(A/\mathfrak{p}_{\nu})$, where the $\mathfrak{p_{\nu}}$ denote the minimal
  prime ideals of $A$, as above. The projection $E\longrightarrow
  E(\boldsymbol{k})$ induces a 
  morphism of functors $\Hom_{A}(-,E) \longrightarrow
  \Hom_{A}(-,E(\boldsymbol{k}))$. This is an isomorphism on modules $M$ of
  finite length, because, for such $M$, $\Hom_{A}(M,A/\mathfrak{p}_{\nu}) =0$. 
\end{proof}

\begin{lemma}\label{completion}
  Let $(A,\mathfrak{m},\boldsymbol{k})$ be a local Noetherian ring and $M$ an
  $A$-module of finite length. Then, there are isomorphisms
  $$M\cong \widehat{M}   \quad \text{ and }\quad
    \mathbb{M}_{A}(M) \cong \mathbb{M}_{\widehat{A}}(\widehat{M})$$
\end{lemma}

\begin{proof}
  We know $A\rightarrow \widehat{A}$ is flat and $\widehat{A}$ is a local
  Noetherian ring with maximal ideal $\mathfrak{m}\widehat{A}$ and residue
  field $\widehat{A}/\mathfrak{m}\widehat{A} \cong \boldsymbol{k}$. In
  particular, $\boldsymbol{k}\otimes_{A}\widehat{A} \cong
  \boldsymbol{k}$. Using these facts, induction on the length $\ell(M)$ shows
  easily that the canonical map $M\rightarrow \widehat{M}\cong
  M\otimes_{A}\widehat{A}$ is an isomorphism.

  From \cite{Matsumura}, Thm. 18.6, we know that the canonical map
  $E(\boldsymbol{k}) \rightarrow E(\boldsymbol{k})\otimes_{A}\widehat{A}$ is an
  isomorphism and that $E(\boldsymbol{k})$ is the injective hull of
  $\boldsymbol{k}$ as an $\widehat{A}$-module as well. This implies that
  $\Hom_{A}(M,E(\boldsymbol{k}))$ and
  $\Hom_{\widehat{A}}(\widehat{M},E(\boldsymbol{k}))$ are isomorphic, which
  gives the claim.
\end{proof}

\begin{remark}\label{reverse}
  Let $(A,\mathfrak{m},\boldsymbol{k})$ be a local Gorenstein ring of
  dimension one, which contains $\boldsymbol{k}$. If $M$ is an $A$-module, we
  can consider it as a $\boldsymbol{k}$-module and equip
  $\Hom_{\boldsymbol{k}}(M,\boldsymbol{k})$ with the structure of an $A$-module
  by $(a\cdot f)(m) = f(am)$. There is a functorial isomorphism of $A$-modules
  \begin{equation}
    \label{homas}
    \Hom_{A}(M,\Hom_{\boldsymbol{k}}(A,\boldsymbol{k})) \xrightarrow{\sim}
    \Hom_{\boldsymbol{k}}(M,\boldsymbol{k}),
  \end{equation}
  which sends $\varphi: M\rightarrow \Hom_{\boldsymbol{k}}(A,\boldsymbol{k})$
  to the mapping $\Phi:M\rightarrow \boldsymbol{k}$ which is defined by
  $\Phi(m):= \varphi(m)(1)$. Because $\Hom_{\boldsymbol{k}}(-,\boldsymbol{k})$
  is exact on the category of $\boldsymbol{k}$-modules it is exact on the
  category of $A$-modules as well, hence (\ref{homas}) implies that the
  $A$-module $\Hom_{\boldsymbol{k}}(A,\boldsymbol{k})$ is injective. If we
  substitute $M=\boldsymbol{k}$ in (\ref{homas}), we obtain
  $$\Hom_{A}(\boldsymbol{k},\Hom_{\boldsymbol{k}}(A,\boldsymbol{k})) \cong
   \Hom_{\boldsymbol{k}}(\boldsymbol{k},\boldsymbol{k}) \cong \boldsymbol{k}.$$
  Now, we can apply lemma \ref{kdual} to
  $E=\Hom_{\boldsymbol{k}}(A,\boldsymbol{k})$ and obtain
  $$\mathbb{M}_{A}(-)\cong \Hom_{\boldsymbol{k}}(-,\boldsymbol{k})$$
  on the category of $A$-modules of finite length. We shall use this
  functorial isomorphism in the case $A=R:=\local{x}{y}$. 

  This isomorphism allows a very nice description of Matlis duality on
  $R$-modules of finite length in terms of the band and string diagrams: all
  the vector spaces are to be replaced by its dual and the mappings by its
  transposed. Because the transposed of $J_{m}(\lambda)$ is similar to
  $J_{m}(\lambda)$, this means: the diagram of the Matlis dual is obtained by
  reversing the direction of all arrows.
\end{remark}

We are interested in a global version of the Matlis functor on a curve. As
usual, if $X$ is a scheme, we denote by $\mathcal{K}$ its sheaf of total
quotient rings. Its stalk $\mathcal{K}_{x}$ is the total quotient ring of the
local ring $\mathcal{O}_{X,x}$.

\begin{lemma}\label{global}
  Let $X$ be a Noetherian scheme of dimension $1$. 
  \begin{itemize}
  \item[(i)] For any closed point $x\in X$ and any $\mathcal{O}_{X,x}$-module
    $W$ we denote by $i_{x}(W)$  the $\mathcal{O}_{X}$-module with stalk $W$ at
    $x$ and $0$ elsewhere. There exists an exact sequence of
    $\mathcal{O}_{X}$-modules
    $$\begin{CD}
      0 @>>> \mathcal{O}_{X} @>>> \mathcal{K} @>>>
      \bigoplus_{x\in X} i_{x}(\mathcal{K}_{x}/\mathcal{O}_{X,x}) @>>> 0,
    \end{CD}$$
    where the direct sum ranges over all closed points $x\in X$. 
  \item[(ii)] If $X$ is Gorenstein, there is an isomorphism of
    $\mathcal{O}_{X}$-modules $$\mathcal{K}/\mathcal{O}_{X} \cong
    \bigoplus_{x\in X} i_{x}E(\boldsymbol{k}(x)).$$
    The direct sum is taken over closed points of $X$, $E(-)$ denotes the
    injective hull of an $\mathcal{O}_{X,x}$-module and
    $\boldsymbol{k}(x)\cong \mathcal{O}_{X,x}/\mathfrak{m}_{X,x}$. 
  \end{itemize}
\end{lemma}

\begin{proof}
  Part (i) is well-known and is obtained by taking stalks. Part (ii) follows
  from (i) and lemma \ref{injective}.
\end{proof}

\begin{definition}
  If $X$ is a Gorenstein curve, we call $$\mathbb{M}(-) :=
  \mathcal{H}om_{X}(-,\mathcal{K}/\mathcal{O}_{X})$$ the (global) \emph{Matlis
  functor}. 
\end{definition}

\begin{remark}
  Let $F$ be a coherent torsion sheaf on a Gorenstein curve $X$ and
  suppose $\supp(F)=\{x\}$. From lemma \ref{global} we obtain: $\mathbb{M}(F)$
  is a sky-scraper sheaf with stalk $\mathbb{M}_{\mathcal{O}_{X,x}}(F_{x})$ at
  $x$. More generally, for any coherent torsion sheaf $F$ and any $x\in X$ one
  has: $$\mathbb{M}(F)_{x} \cong \mathbb{M}_{\mathcal{O}_{X,x}}(F_{x}).$$
  From the results quoted above, we obtain that the canonical map $F
  \rightarrow \mathbb{M}(\mathbb{M}(F))$ is an isomorphism if $F$ is a
  coherent torsion sheaf. This is not true for general sheaves.
  From lemma \ref{injective} we obtain for any torsion sheaf $F$ on $X$:
  $$\mathbb{M}(F) \cong \mathcal{E}xt_{X}^{1}(F,\mathcal{O}_{X}).$$
\end{remark}

\begin{lemma}\label{twist}
  Let $E$ be a semi-stable torsion free sheaf of degree zero on
  $\boldsymbol{E}$ (a Weierstra{\ss} cubic) and define
  $F:=\mathbb{F}(E)$. This is a coherent torsion sheaf on
  $\boldsymbol{E}$. With notation as in section \ref{sec:fm}, we have
  $$T_{\mathcal{O}}(\mathbb{M}(F)) \cong E^{\vee}(-p_{0})[1].$$ 
\end{lemma}

\begin{proof}
  By theorem \ref{ss}, $F$ sits in an exact sequence
  \begin{equation}\label{exsequ}
  \begin{CD}
    0 @>>> H^{0}(E(p_{0}))\otimes\mathcal{O} @>\text{ev}>>
    E(p_{0}) @>>> F @>>> 0.
  \end{CD}
  \end{equation}
  Since $F$ is a torsion sheaf, $\mathbb{M}(F) \cong
  \mathcal{E}xt_{X}^{1}(F,\mathcal{O}_{X})$ is a torsion sheaf as well. In
  particular, $H^{1}(\mathbb{M}(F))=0$ and the evaluation map
  $H^{0}(\mathbb{M}(F))\otimes\mathcal{O} \rightarrow \mathbb{M}(F)$ is
  surjective. Hence, by lemma \ref{computeT}, we have:
  $$T_{\mathcal{O}}(\mathbb{M}(F)) \cong
  \ker(H^{0}(\mathbb{M}(F))\otimes\mathcal{O} \rightarrow \mathbb{M}(F))[1].$$
  Observe $\mathcal{H}om(F,\mathcal{O})=0$ and
  $\mathcal{E}xt^{1}(E(p_{0}),\mathcal{O})=0$. The latter is true because
  $\Ext^{1}_{A}(M,A)=0$ if $A$ is a local Gorenstein ring of dimension $1$ and
  $M$ is a finite $A$-module with $\depth(M)=1$. This follows easily from
  Nakayama's lemma. Hence, application of the functor
  $\mathcal{H}om(-,\mathcal{O})$ to the exact sequence (\ref{exsequ}) gives the
  exact sequence
  $$
    0 \rightarrow E^{\vee}(-p_{0}) \rightarrow
    \mathcal{H}om(H^{0}(E(p_{0}))\otimes\mathcal{O}, \mathcal{O}) 
    \xrightarrow{\delta}
    \mathcal{E}xt^{1}(F,\mathcal{O}) \rightarrow 0.
  $$
  Application of $H^{0}(-)\otimes\mathcal{O}$ gives the following commutative
  diagram:
  $$\begin{CD}
    \mathcal{H}om(H^{0}(E(p_{0}))\otimes\mathcal{O}, \mathcal{O}) @>{\delta}>>
    \mathcal{E}xt^{1}(F,\mathcal{O})\\
    @AA\text{ev}A  @AA\text{ev}A\\
    \Hom(H^{0}(E(p_{0}))\otimes\mathcal{O}, \mathcal{O})\otimes\mathcal{O}
    @>{H^{0}(\delta)\otimes\Id_{\mathcal{O}}}>>
    H^{0}(\mathcal{E}xt^{1}(F,\mathcal{O})) \otimes\mathcal{O}.
  \end{CD}$$
  Obviously, the left vertical evaluation map is an isomorphism.
  Furthermore, the exact sequence
  $$
    H^{0}(E^{\vee}(-p_{0})) 
    \longrightarrow \Hom(H^{0}(E(p_{0})\otimes\mathcal{O}, \mathcal{O})) 
    \xrightarrow{H^{0}(\delta)} H^{0}(\mathcal{E}xt^{1}(F,\mathcal{O}))
  $$
  shows that $H^{0}(\delta)$ is injective, because $\deg(E(p_{0}))>0$ and
  semi-stability of $E(p_{0})$ imply $H^{0}(E^{\vee}(-p_{0})) \cong
  \Hom(E(p_{0}),\mathcal{O}) = 0$.
  As seen before, $\dim \Hom(H^{0}(E(p_{0}))\otimes\mathcal{O},\mathcal{O}) =
  h^{0}(E(p_{0})) = \rk(E)$ as well as $h^{0}(\mathcal{E}xt^{1}(F,\mathcal{O}))
  = \dim \Ext^{1}(F,\mathcal{O}) = \dim \Hom(\mathcal{O},F) = h^{0}(F) =
  \rk(E)$. Hence, $H^{0}(\delta)$ is an isomorphism.
  Now, the conclusion is that $\delta$ and the right vertical evaluation map in
  the diagram above have isomorphic kernels, that is:
  $E^{\vee}(-p_{0})[1] \cong T_{\mathcal{O}}(\mathbb{M}(F)).$
\end{proof}

\begin{definition}\label{twistedmd}
  If $F$ is a torsion sheaf on $\boldsymbol{E}$, we define:
  $$F^{\star} := i^{\ast}\mathbb{M}(F)$$ and call it the \emph{twisted Matlis
    dual} of $F$.
\end{definition}

\begin{theorem}\label{dual}
  If $E$ is a semi-stable torsion free sheaf of degree zero on
  $\boldsymbol{E}$, there is an isomorphism:
  $$\mathbb{F}(E^{\vee}) \cong \mathbb{F}(E)^{\star}.$$
\end{theorem}

\begin{proof}
  Using lemma \ref{twist}, we obtain $E^{\vee}[1]\cong
  T_{\boldsymbol{k}(p_{0})} T_{\mathcal{O}}(\mathbb{M}(F))$. As
  $F=\mathbb{F}(E)$ and $\mathbb{M}(F)$ are torsion sheaves,
  $T_{\boldsymbol{k}(p_{0})} \mathbb{M}(F) \cong \mathbb{M}(F)$ and we obtain
  $E^{\vee}[1]\cong \mathbb{F}\mathbb{M}(F) \cong
  \mathbb{F}\mathbb{M}\mathbb{F}(E)$. This implies $\mathbb{F}(E^{\vee})[1]
  \cong \mathbb{F}\mathbb{F}\mathbb{M}\mathbb{F}(E)$. From theorem \ref{main}
  we deduce $\mathbb{F}(E^{\vee}) \cong i^{\ast}\mathbb{M}(\mathbb{F}(E))
  \cong  \mathbb{F}(E)^{\star}.$
\end{proof}

\begin{remark}
  If $\boldsymbol{E}$ is a smooth elliptic curve, our result is a special case
  of a result of Mukai \cite{Mukai}, (3.8). 
\end{remark}

\begin{remark}
  As before, $R$ denotes the complete local Gorenstein ring $\local{x}{y}$. If
  the singular point $s\in\boldsymbol{E}$ is a node, the ring map
  $R\rightarrow R$, which is induced by $i:\boldsymbol{E}\rightarrow
  \boldsymbol{E}$ at the fixed point $s\in\boldsymbol{E}$, interchanges $x$
  and $y$. Hence, if $M$ is an $R$-module of finite length, the diagram for
  $i^{\ast}(M)$ is obtained from the diagram of $M$ by interchanging $x$ and
  $y$. Let now $E$ be an indecomposable semi-stable torsion free sheaf of
  degree zero. Using lemma \ref{compDual} and theorem \ref{description} it is
  not hard to verify that the diagram of $i^{\ast}\mathbb{F}(E^{\vee})$ is
  obtained from that of $\mathbb{F}(E)$ by reversing all arrows. Because, by
  theorem \ref{dual}, $\mathbb{M}_{R}(\mathbb{F}(E)) 
  \cong i^{\ast}\mathbb{F}(E^{\vee})$, this gives another verification of the
  description of the Matlis dual of an $R$-module of finite length which was
  given in remark \ref{reverse}. We illustrate this by the two examples below.
\end{remark}

\begin{example}
  The Matlis dual module of the string module $$\mathcal{N}(2(3,2)1) \cong
  \mathbb{F}(\mathcal{S}(-1,0,1,0,0,-1,0,1,-1))$$ is the module
  \begin{align*}
   i^{\ast}\mathbb{F}(\mathcal{S}(-1,0,1,0,0,-1,0,1,-1)^{\vee}) &\cong
   i^{\ast}\mathbb{F}(\mathcal{S}(0,0,-1,0,0,1,0,-1,0))   \\ \cong
   \mathbb{F}(\mathcal{S}(0,-1,0,1,0,0,-1,0,0))     &\cong
   \mathcal{N}(0(1,2)(3,2)0)
  \end{align*} This module is represented by the diagram which is obtained
  from the diagram in example \ref{ex:string} by reversing all arrows:

  
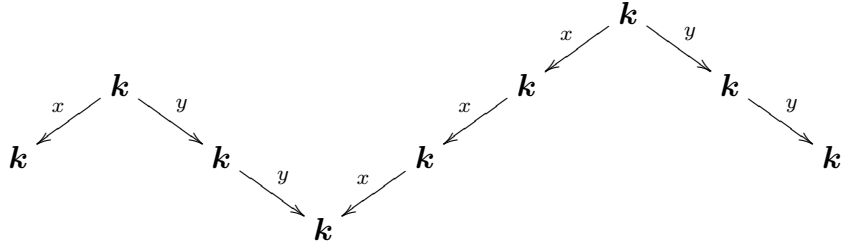
\xymatrix@R-=12pt{
&&&&&&{\boldsymbol{k}}\ar@{->}[dl]_{x}\ar@{->}[dr]^{y}\\
& {\boldsymbol{k}}\ar@{->}[dl]_{x}\ar@{->}[dr]^{y}&&&&   
  {\boldsymbol{k}}\ar@{->}[dl]_{x} && 
  {\boldsymbol{k}}\ar@{->}[dr]^{y} \\
  {\boldsymbol{k}}&&   {\boldsymbol{k}}\ar@{->}[dr]^{y}&& 
  {\boldsymbol{k}}\ar@{->}[dl]_{x}&&&& {\boldsymbol{k}} \\
&&&{\boldsymbol{k}}  
}
\end{example}

\begin{example}
  The Matlis dual module of the band module
  \begin{align*}
   &\mathcal{M}(((2,2)(3,4)(1,3)),m,\lambda) \\
   \cong &\mathbb{F}(\mathcal{B}((1,0,-1,0,1,0,0,-1,0,0,0,1,-1,0,0),m,\lambda))
  \end{align*}
  is the module
  \begin{align*}
    &i^{\ast}
   \mathbb{F}(\mathcal{B}((1,0,-1,0,1,0,0,-1,0,0,0,1,-1,0,0),m,\lambda)^{\vee})
     \\
    \cong&i^{\ast}
   \mathbb{F}(\mathcal{B}((-1,0,1,0,-1,0,0,1,0,0,0,-1,1,0,0),m,\lambda^{-1}))
     \\
    \cong&
   \mathbb{F}(\mathcal{B}((0,0,1,-1,0,0,0,1,0,0,-1,0,1,0,-1),m,\lambda))\\
    \cong&
   \mathbb{F}(\mathcal{B}((1,-1,0,0,0,1,0,0,-1,0,1,0,-1,0,0),m,\lambda))\\
    \cong&\mathcal{M}((1,4)(3,2)(2,3)),m,\lambda)
  \end{align*}
  This module is represented by the following diagram which is obtained
  from the diagram in example \ref{ex:band} by reversing all arrows (and
  moving the Jordan block $J_{m}(\lambda)$ two arrows forward): 

\xymatrix@C-=4pt{
&{\boldsymbol{k}^{m}}\ar@{->}[dl]_{xI_{m}}\ar@{->}[dr]_{yI_{m}}\\
{\boldsymbol{k}^{m}}
&&{\boldsymbol{k}^{m}}\ar@{->}[dr]_{yI_{m}}
&&&&&&{\boldsymbol{k}^{m}}\ar@{->}[dl]_{xI_{m}}\ar@{->}[dr]_{yI_{m}}
&&&&{\boldsymbol{k}^{m}}\ar@{->}[dl]_{xI_{m}}\ar@{->}[d]^{yI_{m}}\\
&&&{\boldsymbol{k}^{m}}\ar@{->}[dr]_{yI_{m}}
&&&&{\boldsymbol{k}^{m}}\ar@{->}[dl]_{xI_{m}}
&&{\boldsymbol{k}^{m}}\ar@{->}[dr]_{yI_{m}}
&&{\boldsymbol{k}^{m}}\ar@{->}[dl]_{xI_{m}}
&{\boldsymbol{k}^{m}}\ar@{->}[d]^{yI_{m}}\\
&&&&{\boldsymbol{k}^{m}}\ar@{->}[dr]_{yI_{m}}
&&{\boldsymbol{k}^{m}}\ar@{->}[dl]_{xI_{m}}
&&&&{\boldsymbol{k}^{m}}
&&{\boldsymbol{k}^{m}}\ar@/^8pc/@{->}[uullllllllllll]_(.3){yJ_{m}(\lambda)}\\
&&&&&{\boldsymbol{k}^{m}}
}
\end{example}

\end{document}